\documentclass{amsart} 

\usepackage{amsmath,amsthm}     
\usepackage{amssymb}            
\usepackage{euscript}           
\usepackage{enumerate,calc}
\usepackage[matrix,arrow,curve,frame]{xy}    

\xymatrixcolsep{1.9pc}                          
\xymatrixrowsep{1.9pc}
\newdir{ >}{{}*!/-5pt/\dir{>}}                  

\newtheorem{thm}[subsection]{Theorem}
\newtheorem{defn}[subsection]{Definition}
\newtheorem{prop}[subsection]{Proposition}

\newtheorem{cor}[subsection]{Corollary}
\newtheorem{lemma}[subsection]{Lemma}

\theoremstyle{definition}  
\newtheorem{ex}[subsection]{Example}
\newtheorem{exercise}[subsection]{Exercise}

\newtheorem{remark}[subsection]{Remark}

  {\end{list}}
  
\newcommand{\dfn}{\textbf} 

\newcommand{\mdfn}[1]{\dfn{\mathversion{bold}#1}} 

\newcommand{\field}[1]  {\mathbb #1} 

\newcommand{\R}         {\field R}

\newcommand{\HH}        {\field H}

\newcommand{\C}         {\field C}
\renewcommand{\O}         {\field O}

\DeclareMathOperator{\Ann}{Ann}
\DeclareMathOperator{\Alt}{Alt}
\DeclareMathOperator{\Ass}{Ass}

\newcommand{\ra}{\rightarrow}                   

\newcommand{\rea}[1]{|{#1}|}             
\newcommand{\map}{\rightarrow}

\newcommand{\ceck}[1]{\Cech(#1)}         
\newcommand{\oceck}[1]{\Cech^{o}(#1)}    
\newcommand{\oreal}[1]{\rea{\oceck{U}}}  
\newcommand{\creal}[1]{\rea{\ceck{U}}}   

\newcommand{\Cech}{\check{C}}

\newcommand{\ZD}{ZD}

\renewcommand{\Re}{\text{Re}}
\renewcommand{\Im}{\text{Im}}

\numberwithin{equation}{subsection}

\newenvironment{myequation}
  {\addtocounter{subsection}{1}\begin{eqnarray}}
  {\end{eqnarray}$\!\!$}

\begin{document}

\title{Large annihilators in Cayley-Dickson algebras}

\author{Daniel K. Biss}\thanks{This research was conducted during the
period the first author served as a Clay Mathematics Institute
Research Fellow}
\author{Daniel Dugger}
\author{Daniel C. Isaksen}

\address{Department of Mathematics\\
University of Chicago\\Chicago, IL 60637}

\address{Department of Mathematics\\ University of Oregon\\ Eugene, OR
97403}

\address{Department of Mathematics\\ Wayne State University\\
Detroit, MI 48202}

\email{daniel@math.uchicago.edu}

\email{ddugger@math.uoregon.edu}

\email{isaksen@math.wayne.edu}

\begin{abstract}
Cayley-Dickson algebras are non-associative $\R$-algebras that generalize
the well-known algebras $\R$, $\C$, $\HH$, and $\O$.  We study zero-divisors
in these algebras.  In particular, we show that the annihilator of any element
of the $2^n$-dimensional Cayley-Dickson algebra has dimension at most
$2^n - 4n +4$.  Moreover, every multiple of $4$ between $0$ and this
upper bound occurs as the dimension of some annihilator.  Although a 
complete description of zero-divisors seems to be out of reach, we
can describe precisely the elements whose annihilators have dimension
$2^n - 4n + 4$.
\end{abstract}

\maketitle

\section{Introduction}

The Cayley-Dickson algebras are a sequence $A_0,A_1,\dots$ of
non-associative $\R$-algebras with involution.  The first few are
familiar: $A_0 = \R,$ $A_1 = \C$, $A_2 = \HH$ (the quaternions), and
$A_3 = \O$ (the octonions).  Each algebra $A_n$ is constructed from
the previous one $A_{n-1}$ by a doubling procedure; unfortunately,
this doubling procedure tends to destroy desirable algebra properties.
For example, $\R$ is the only Cayley-Dickson algebra with trivial
involution, $\R$ and $\C$ are the only commutative
Cayley-Dickson algebras, and $\R$, $\C$, and $\HH$ are the only
associative ones.

Only the first four Cayley-Dickson algebras, $\R$, $\C$, $\HH$, and
$\O$, are normed; equivalently, no $A_n$ with $n\geq 4$ is alternative.
When $n=4$ two things happen.  First of all, this weakening
sequence of algebraic conditions stops: there is not an easily-phrased
algebraic condition satisfied by $A_n$ for $n\leq 4$ but not $A_5$,
or $A_n$ for $n\leq 5$ but not $A_6$, and so forth.  Secondly, for all
$n\geq 4$, the algebra $A_n$ admits non-trivial zero-divisors.  That
is, there exist pairs of non-zero elements $x$ and $y$ in $A_n$ such that 
$xy = 0$.

Historically, these two mysterious facts have conspired to discourage
mathematicians from studying the higher Cayley-Dickson algebras.
We instead take them as
our point of departure:  the locus of zero-divisors in $A_4$ is
well-understood and both quite simple and interesting.  In our view,
the zero-divisors in the higher $A_n$ offer some promise in two
distinct ways.  First, as $n$ grows the locus of zero-divisors in
$A_n$ becomes more complicated, and this complexity serves as an
analogue to the weakening sequence of algebraic criteria present in
$A_0,A_1, A_2$, and $A_3$.  Secondly, it is our hope that these loci
of zero-divisors will prove to be geometrically interesting in their own right.

Accordingly, the goal of this article is to initiate a systematic
study of the zero-divisors in the Cayley-Dickson algebras.

\medskip

\subsection{Statement of results}
We now state the main results in more detail.  If $x$ belongs to $A_n$,
then the
\dfn{annihilator} of $x$ is $\Ann x=\{y\in A_n : xy=0\}$.  Because $A_n$ is
noncommutative one should really call this the \dfn{right annihilator}
of $x$, but it turns out that in Cayley-Dickson algebras the left and
right annihilators of an element are always the same (see
Corollary~\ref{cor:left=right}).  

A theorem of Moreno \cite{M1} says that the (real) dimension of $\Ann
x$ is always a multiple of $4$.  The reader will find a different
proof of this result in our Theorem~\ref{thm:4-dim}.  The first goal
of this paper is to determine exactly which multiples of $4$ can
occur:

\begin{thm} 
If $x$ belongs to $A_n$,
then $\dim (\Ann x) \leq 2^n-4n+4$.  Moreover, if $d$ is
any multiple of $4$ such that $0\leq d \leq 2^n-4n+4$, then there exist
elements of $A_n$ whose annihilators have dimension $d$.
\end{thm}

So in $A_4$ there are annihilators of dimensions $0$ and $4$ only; in
$A_5$ there are annihilators of dimensions $0, 4, 8, 12,$ and $16$; and 
in $A_6$ there are annihilators of dimensions $0, 4, 8, \ldots,
44$.  Notice that the maximal dimension of the annihilators grows very
quickly, the codimension in $A_n$ being given by a linear function.
When $n$ is large, one has zero-divisors whose annihilator is
`almost' the whole algebra.

Also note that the above result even makes sense when $n\leq 3$, in
which case it gives the well-known fact that these algebras have no 
zero-divisors.

\medskip

The space of zero-divisors in $A_n$ is closed under scalar
multiplication, and so it forms a cone in the real vector space
underlying $A_n$.  It is therefore natural to focus on the norm $1$
elements and look at the space
\[
\ZD(A_n)=\{x\in A_n :  ||x||=1, \ \Ann x \neq 0\}.
\]  
We would eventually like to understand the topological properties of
this space, although at the moment this seems to be a very complicated
problem.  It is natural, perhaps, to look at the subspaces
\[ \ZD_k(A_n)=\{x\in A_n :  ||x||=1, \ \dim (\Ann x)=k\}.
\]
These strata are also complicated, and unknown even in the case of
$A_5$.   

The dimension---and complexity---of the strata increases as $k$
becomes small.  So the `simplest' stratum is $\ZD_{2^n-4n+4}(A_n)$.
With some effort we can describe this one completely:

\begin{thm}
\label{thm:topZD}
When $n\geq 4$ the space $\ZD_{2^n-4n+4}(A_n)$ is homeomorphic to a
disjoint union of $2^{n-4}$ copies of the Stiefel variety $V_2(\R^7)$,
i.e., the space of ordered pairs of orthonormal vectors in $\R^7$.
\end{thm} 

Obtaining a description of $\ZD_k(A_n)$ for smaller $k$ is an
intriguing open problem.

\medskip

We close this introduction by describing one technique which is used
repeatedly in the paper, and which will help explain
Theorem~\ref{thm:topZD}.  Every Cayley-Dickson algebra $A_n$ contains
a distinguished element $i_n$, and the $2$-dimensional subspace
$\langle 1,i_n\rangle$ is a subalgebra isomorphic to $\C$ that we
denote by $\C_n$.
It turns out that $A_n$ is a vector space over $\C_n$.
The special properties of $i_n$ guarantee that multiplication
in $A_n$ behaves well with respect to $\C_n$; see 
Lemmas \ref{lem:i-comm} and \ref{lem:i-comm2} for precise statements.
This viewpoint on $A_n$ underlies almost all of our main results.

For every $a$ in $A_n$ that is orthogonal to $\C_n$, we show in
Theorem~\ref{thm:C-ann2} that $(a,\pm i_na)$ is a zero-divisor in
$A_{n+1}$.  Moreover, we completely determine its annihilator and find
that it has dimension $2^n-4+\dim (\Ann a)$.  So this gives a method
for producing `large' annihilators inside of $A_{n+1}$.

The space $\ZD_4(A_4)$ is easily analyzed by hand and shown to be $V_2(\R^7)$.
Applying the maps $a\mapsto (a,\pm i_4a)$ gives two disjoint copies of
$V_2(\R^7)$ inside of $\ZD_{16}(A_5)$, and then applying $a \mapsto (a,\pm
i_5a)$ gives four disjoint copies inside of $\ZD_{44}(A_{6})$.  This
explains the spaces arising in Theorem~\ref{thm:topZD}. The proof of
that theorem involves showing that when $n\geq 5$ {\it every\/}
$(2^n-4n+4)$-dimensional annihilator in $A_n$ is of the form 
$\Ann(a,\pm i_{n-1}a)$ for some $a$ in $A_{n-1}$.

\subsection{Presentation of the paper}
It is our intention that this paper be the first in a series, and so
we have been extremely careful in laying out the foundations.  For
this reason we have occasionally chosen to duplicate known results,
either because they are not published or because the published proofs
do not fit well into our approach.  In particular, Guillermo Moreno has
a number of nice results on zero-divisors in Cayley-Dickson
algebras~\cite{M1}, some of which we have chosen to reprove here in the
interests of building an arsenal of techniques.  We have made every
effort to state clearly which theorems were already known to Moreno.

One benefit is that the text should be completely self-contained,
except for undergraduate-level linear algebra and abstract algebra,
and a small amount of topology.

The other guiding principle behind our expository style has been to
avoid ad hoc proofs and results.  We'll try to use the same basic
ideas in all of the proofs.  When a phenomenon occurs, we'll emphasize
the underlying structure that makes it happen, rather than just record
the technical consequences.  We have not completely attained this
goal, however.
For instance, it breaks down in our proof of
Theorem~\ref{thm:mainA_5}, where we arbitrarily manipulate
linear quaternionic equations.

There are a variety of other papers in the literature on
Cayley-Dickson algebras.  Dickson's contribution was to construct the
octonions (also sometimes called the Cayley numbers) as pairs of
quaternions \cite{Di}.  It was actually Adrian Albert (a student of
Dickson) who first iterated this construction to produce an infinite
family of algebras \cite{Al}.  Since then, a number of articles have
been written studying various algebraic properties of these algebras.
We mention in particular \cite{A}, \cite{Br}, \cite{ES}, \cite{KY},
\cite{M1}, and \cite{Sc}.  Fred Cohen has suggested that an
understanding of the zero-divisor loci of the Cayley-Dickson algebras
might have useful applications in topology---see \cite{Co} for more
information.

\subsection{Acknowledgments} 
The authors would like to thank Dan Christensen for his assistance with
several computer calculations.

\section{Cayley-Dickson algebras}

We start with the inductive definition of Cayley-Dickson algebras.
These are finite-dimensional $\R$-algebras equipped with a linear 
involution $(-)^*$ satisfying $(xy)^*=y^*x^*$.  The few results of this
section appear to have been known to every mathematician from the
modern era who worked on Cayley-Dickson algebras.

\begin{defn}
\label{defn:CD}
The algebra \mdfn{$A_0$} is equal to $\R$, and the involution is the identity.
The algebra \mdfn{$A_n$} additively is 
$A_{n-1} \times A_{n-1}$.  The multiplication is defined inductively by
\[
(a,b)(c,d) = (ac - d^* b, da + bc^*),
\]
and the involution is also defined inductively by
\[
(a,b)^* = (a^*, -b).
\]
\end{defn}

Note that $A_n$ has dimension $2^n$ as an $\R$-vector space.  We shall
use the term ``conjugation'' to refer to the involution because it
generalizes the usual conjugation on the complex numbers.  The reader
should verify inductively that conjugation is in fact an
involution---that is, \mdfn{$(x^*)^* = x$}---and that it interacts with
multiplication according to the formula \mdfn{$(xy)^* = y^* x^*$}.

One can see inductively that $A_n$ contains a copy of $\R$
because the subalgebra $A_{n-1} \times 0$ of $A_n$ is isomorphic
to $A_{n-1}$.  Moreover, since the inductive formulas defining multiplication
and conjugation are $\R$-linear, we see that each $A_n$ is an $\R$-algebra.
In particular, the elements of $\R$ are central in $A_n$.

\begin{ex}
The algebra $A_1$ is isomorphic to the complex numbers $\C$ with
its usual conjugation.  To see why, just check that $(0,1)$ plays the
role of $i$ in $\C$.
\end{ex}

\begin{ex}
The algebra $A_2$ is isomorphic to the quaternions $\HH$ with its usual
conjugation.  The elements $(i,0)$, $(0,1)$, and $(0,i)$ play
the roles of the standard basis elements $i$, $j$, and $k$.
\end{ex}

\begin{ex}
The algebra $A_3$ is isomorphic to the octonions $\O$ with its usual
conjugation.  
\end{ex}

\begin{defn}
If $x$ is any vector, then the \mdfn{real part $\Re(x)$} of $x$
is defined inductively as follows.  If $x$ belongs to $A_0$, then $\Re(x)$
equals $x$.  If $(a,b)$ belongs to $A_n$, then $\Re(a,b) = \Re(a)$.
Also, the \mdfn{imaginary part $\Im(x)$} of $x$ 
is equal to $x - \Re(x)$.
\end{defn}

If $x$ belongs to $\R$,
then we say that $x$ is \mdfn{real}.  If $\Re(x) = 0$, then we say that
$x$ is \mdfn{imaginary}.  Every vector can be uniquely written as the sum
of a real vector and an imaginary vector.

The reader should check inductively that \mdfn{$2\Re(x) = x + x^*$} and
\mdfn{$2\Im(x) = x - x^*$}.
One consequence is that \mdfn{$x^* = -x$ if $x$ is imaginary}.

\begin{lemma}
\label{lem:Re-comm}
For all $x$ and $y$, $\Re(xy - yx) = 0$.
\end{lemma}

\begin{proof}
If either $x$ or $y$ is real, then the formula is trivial.  By linearity,
we may assume that $x$ and $y$ are both imaginary, so $x^* = -x$
and $y^* = -y$.  Therefore,
\[
xy - yx = xy - y^* x^* = xy - (xy)^* = 2\Im(xy).
\]
Thus, $xy - yx$ is imaginary, so its real part vanishes.
\end{proof}

\begin{defn}
\label{defn:assoc}
The \mdfn{associator} of $x$, $y$, and $z$ is
\mdfn{$[x,y,z]= (xy)z - x(yz)$}.
\end{defn}

\begin{lemma}
\label{lem:Re-assoc}
For all $x$, $y$, and $z$, $\Re([x,y,z]) = 0$.
\end{lemma}

\begin{proof}
By linearity, it suffices to assume that $x$, $y$, and $z$ are of
the form $(a,0)$ or $(0,a)$.  There are eight cases to check.
Four of the cases are easy because both $(xy)z$ and $x(yz)$ are of
the form $(0,a)$; thus their real parts are zero.

The other four cases have to be checked one at a time.  The case 
$x = (a,0)$, $y = (b,0)$, $z = (c,0)$ follows by induction.  The other
three cases are all very similar to one another; we give one example.

Let $x = (a,0)$, $y = (0,b)$, and $z = (0,c)$.  Then calculate that
$[x,y,z] = ( a(c^* b) - c^*(ba), 0)$.
Now 
$\Re( a(c^* b) - c^*(ba) )$ equals $\Re( (c^* b) a - c^*(ba))$
by Lemma~\ref{lem:Re-comm}, and this last expression
vanishes by induction.
\end{proof}

One might be tempted to conclude from Lemmas~\ref{lem:Re-comm} and
\ref{lem:Re-assoc} that when computing the real part of any expression,
one can completely ignore the order of the factors and the parentheses.
This is not true.  For example, we cannot conclude that $\Re(x(yz))$
equals $\Re(x(zy))$.

The following corollary says that Cayley-Dickson algebras are ``flexible''.
This means that expressions of the form $xyx$ are well-defined, even
in the absence of associativity.

\begin{cor}
\label{cor:flex}
For all $x$ and $y$ in $A_n$, $[x,y,x]$ vanishes.
\end{cor}

\begin{proof}
It suffices to prove the identity
\[
x(yz) - (xy)z + z(yx) - (zy)x = 0
\]
because we can set $z = x/2$ to recover the original formula.

If any one of $x$, $y$, and $z$ is real, then the formula is trivial.
By linearity, we may assume that $x$, $y$, and $z$ are imaginary.
Using that $x^* = -x$, $y^* = -y$, and $z^* = -z$, 
the above expression is its own conjugate.
Therefore, the imaginary part of the expression
is zero.

On the other hand, the real part of the expression is also zero
by Lemma~\ref{lem:Re-assoc}.
\end{proof}

\begin{defn}
\label{defn:standard-basis}
The \mdfn{standard basis} of $A_n$ is defined inductively as follows.
The standard basis for $A_0$ consists of the element $1$.  The standard
basis of $A_n$ consists of elements of the form $(x,0)$ or $(0,x)$,
where $x$ belongs to the standard basis of $A_{n-1}$.
\end{defn}

One readily checks by induction that the above does
indeed give an $\R$-vector space basis of $A_n$.

\section{Real inner product}

The definition and basic properties of the real inner product
described in this section were laid out in \cite{M1}.

\begin{defn}
\label{defn:inner-product}
The \mdfn{real inner product $\langle x, y \rangle$} 
of two elements $x$ and $y$ in $A_n$ is equal to $\Re(xy^*)$.
\end{defn}

\begin{prop}
\label{prop:inner-product}
The function $\langle -,- \rangle$ is a symmetric positive-definite inner
product on $A_n$.
\end{prop}

\begin{proof}
The function $\langle -,- \rangle$ is $\R$-bilinear.  For symmetry,
$\Re(xy^*)$ equals $\Re(yx^*)$ because
conjugation does not change the real part of a vector.

For positive-definiteness,
compute that
\[
\Re((a,b)(a,b)^*) = \Re(aa^* + b^* b).
\]
By induction, $\Re(aa^*) + \Re(bb^*) \geq 0$, and 
$\Re(aa^*) + \Re(bb^*) = 0$ if and only if $a = b = 0$.
\end{proof}

As usual, we set $||x||=\sqrt{\langle
x,x\rangle}=\sqrt{\Re(xx^*)}$.  This makes sense because of
positive-definiteness.  

Recall from Definition~\ref{defn:standard-basis} that we can identify
$A_n$ with $\R^{2^n}$ (as $\R$-vector spaces).  It is easy to check by
induction that under this identification,
the inner product of Definition~\ref{defn:inner-product} corresponds to
the standard inner product on $\R^{2^n}$.  Therefore, the standard
basis of $A_n$ is in fact an orthonormal basis.

\begin{defn}
\label{defn:L_x}
For any vector $x$, let \mdfn{$L_x$} be the linear endomorphism 
of $A_n$ given by left multiplication by $x$, and 
let \mdfn{$R_x$} be the linear endomorphism 
of $A_n$ given by right multiplication by $x$.
\end{defn}

\begin{lemma}
\label{lem:adjoint}
The maps $L_x$ and $L_{x^*}$ are adjoint in the sense that
$\langle L_x y, z \rangle$ equals $\langle y, L_{x^*} z \rangle$.
The maps $R_x$ and $R_{x^*}$ are also adjoint.
\end{lemma}

\begin{proof}
For the first claim, 
we want to show that $\Re((xy)z^*)$ equals $\Re(y(z^* x))$.
This follows from Lemma~\ref{lem:Re-comm}, Lemma~\ref{lem:Re-assoc},
and then Lemma~\ref{lem:Re-comm} again.

The proof for $R_x$ is similar.
\end{proof}

\begin{cor}
\label{cor:L_x-sym}
Suppose that $x$ is imaginary.
The map $L_x$ is antisymmetric in the sense that
$\langle L_x y, z \rangle = - \langle y, L_x z \rangle$.
The map $R_x$ is also antisymmetric.
\end{cor}

\begin{proof}
This follows immediately from Lemma~\ref{lem:adjoint}
since $x^* = -x$.
\end{proof}

Note that $L_x$ need not be an isometry, even if $x$ has norm $1$.
For example, if $x$ is a zero-divisor then there exists a non-zero $y$
such that $xy = 0$.  Then $\langle y, y \rangle$ is not zero, but
$\langle L_x y, L_x y \rangle$ is zero.

\begin{lemma}
\label{lem:norm}
For any $x$, $x x^* = x^* x = ||x||^2$.
\end{lemma}

\begin{proof}
First note that $x x^*$ is real because it is its own conjugate.
Then $x x^*$ equals $\Re(x x^*)$, which by definition is $||x||^2$.
Finally, note that $x^* x$ equals $\Re(x^* x)$; now apply
Lemma~\ref{lem:Re-comm} to show that $x^* x$ equals $x x^*$.
\end{proof}

One consequence is that \mdfn{$x^2 = -||x||^2$ if $x$ is imaginary}.
Thus, the square of an imaginary vector is zero if 
and only if the original vector is zero.

\begin{lemma}
\label{lem:anti-comm}
Two imaginary vectors $x$ and $y$ anti-commute if and only if 
they are orthogonal.
\end{lemma}

\begin{proof}
Since $x^* = -x$ and $y^* = -y$, the conjugate of $xy$ is $yx$.  
Therefore, $x$ and $y$ anti-commute if and only if $xy$ is imaginary,
i.e., $\Re(xy)$ is zero.  Finally, $\Re(xy)$ equals $\langle
x,y^*\rangle = -\langle x, y \rangle$.
\end{proof}

Throughout the text, if $V$ is a vector space with inner product and
$W$ is a subspace, we denote by $W^{\perp}$ the orthogonal complement
of $W$ in $V$.

\begin{lemma} 
\label{lem:subalg}
Suppose $B$ is an $\R$-subalgebra of $A_n$.
If $b$ is in $B$ and $x$ is in $B^\perp$,
then $bx$ and $xb$ are in $B^\perp$.
\end{lemma}

\begin{proof}
If $b$ lies in $B$, then $b^*$ also lies in $B$ 
because $b^* = -b + 2\Re(b)$ and $\Re(b)$ lies in $B$.
We need to show that $\langle a, bx \rangle$ equals zero 
for every $a$ in $B$.  This equals $\langle b^*a, x \rangle$
by Lemma~\ref{lem:adjoint}, which is zero because 
$b^* a$ belongs to $B$.

A similar argument shows that $xb$ is also in $B^\perp$.  
\end{proof}

The algebras $A_0$, $A_1$, $A_2$, and $A_3$ are normed in the sense that
$||xy|| = ||x|| \cdot ||y||$.  However, when $n \geq 4$ the algebra $A_n$
is not normed---the presence of zero-divisors prevents this.    

\begin{exercise}
If $a,x\in A_n$ and $[a,a,x]=0$, prove using Lemmas~\ref{lem:Re-comm}
and ~\ref{lem:Re-assoc} that $\langle ax,ay\rangle =
||a||^2 \langle x,y\rangle$ for any $y\in A_n$.  
Since $[a,a,x]$ always vanishes in $A_3$ (see
Lemma~\ref{lem:A_3-alt}),
this shows that $A_3$ is normed.
\end{exercise}

\section{Alternators}
\label{sctn:alt}

Basic statements of alternativity, including our Lemma
\ref{lem:alt-basis}, were established by Schafer \cite{Sc}.  Moreno
was able to completely classify the alternative elements in every
$A_n$ \cite{M2}.

\begin{defn}
\label{defn:Alt}
For $x$ in $A_n$, let \mdfn{$\Alt x$} be the linear subspace of all $y$
such that $[x,x,y]$ vanishes.
\end{defn}

We call this space the ``alternator'' of $x$.

\begin{defn}
\label{defn:alt}
An element $x$ of $A_n$ is \mdfn{alternative} if 
$\Alt(x)$ equals $A_n$, i.e., if all expressions $[x,x,y]$ vanish.
\end{defn}

In $A_0$, $A_1$, and $A_2$, all elements are alternative; this follows
immediately from the fact that these algebras are associative.  Even
though $A_3$ is not associative, it turns out that every vector in
$A_3$ is still alternative.

\begin{lemma}
\label{lem:A_3-alt}
Every element of $A_3$ is alternative.
\end{lemma}

\begin{proof}
Simply compute that $[(a,b),(a,b),(x,y)] = 0$ for all quaternions
$a$, $b$, $x$, and $y$.  One needs that $a + a^*$ and $bb^*$ are
both real and therefore central.
\end{proof}

\begin{lemma}
\label{lem:alt-basis}
Every standard basis vector (see Definition \ref{defn:standard-basis})
is alternative.
\end{lemma}

\begin{proof}
The proof is by induction.  Direct computation shows that the associators
$[(x,0),(x,0),(a,b)]$ and $[(0,x),(0,x),(a,b)]$ both vanish whenever
$x$ is alternative in $A_{n-1}$.  For the second equation, the calculation
is simplified by writing $x$ as the sum of a real and an imaginary vector.
\end{proof}

\begin{remark}
It is natural to ask what are the possible
dimensions of $\Alt x$, for $x$ in $A_n$.  Although we will not need
this in the present paper, we can show that $\dim(\Alt x)$ is always a
multiple of $4$ (compare Theorem \ref{thm:4-dim}).
In $A_4$, the only two possible dimensions are $8$ and $16$.
In $A_5$, computer calculations show that there are examples of $x$
with $\dim(\Alt x)$ equal to $4$, $8$, $12$, $16$, $24$, and $32$.
We do not know if dimensions $20$ and $28$ can occur.  
\end{remark}

Later we will use the following definition frequently when considering
specific examples.

\begin{defn}
\label{defn:quat-pair}
Two vectors $a$ and $b$ in $A_n$ are a \dfn{quaternionic pair} if 
there is an injective algebra map $\phi:\HH \map A_n$ such that
$\phi(i) = a$ and $\phi(j) = b$.
\end{defn}

\begin{remark}
Vectors $a$ and $b$ are a quaternionic pair if and only if they are
orthogonal imaginary unit vectors such that $[a,b,b]$ and $[a,a,b]$
both vanish.  This is because $\HH$ is the free associative
$\R$-algebra on $i$ and $j$ subject to the relations $i^2 = j^2 = -1$
and $ij = -ji$.  Lemma~\ref{lem:norm} tells us that $a^2$ and $b^2$
are both equal to $-1$ if and only if they are both imaginary unit
vectors, and Lemma~\ref{lem:anti-comm} tells us that the imaginary unit
vectors $a$ and $b$ anti-commute if and only if they are orthogonal.
Finally, in the presence of flexibility of $A_n$, the vanishing of
$[a,a,b]$ and $[b,b,a]$ is equivalent to associativity of the algebra
generated by $a$ and $b$. 
\end{remark}

When $a$ and $b$ are a quaternionic pair,
we write \mdfn{$\HH \langle a,b \rangle$}
 for the subalgebra of $A_n$ that they generate.
Additively, $\HH\langle a,b \rangle$ has an orthonormal basis consisting of 
$1$, $a$, $b$, and $ab$.

\begin{ex}
Let $a$ and $b$ be any two distinct imaginary standard basis elements of $A_n$.
Then $a$ and $b$ are a quaternionic pair.  
See Lemma~\ref{lem:alt-basis}
for the fact that $[a,a,b]$ and $[b,b,a]$ both vanish.
\end{ex}

\section{Complex structure}
\label{sctn:complex}

\begin{defn}
\label{defn:C}
Let \mdfn{$i_n$} be the element $(0,1)$ of $A_n$.
Let \mdfn{$\C_n$} be the subalgebra of $A_n$ additively generated by $1$ 
and $i_n$.
\end{defn}

The notation suggests that $i_n$ is the $n$th analogue
of the square root of $-1$ in $\C$.
Note that $\C_n$ is isomorphic to the complex numbers, 
where $i_n$ plays the role of $i$.  Our first goal is to show that
$A_n$ is a complex vector space, where the $\C$-action is 
given by left multiplication by elements of $\C_n$.

\begin{lemma}
\label{lem:i_n-alt}
For any $x$ in $A_n$, $[i_n,i_n,x]=0$.
\end{lemma}

\begin{proof}
This follows by direct computation with the definition of
multiplication (see also Lemma~\ref{lem:alt-basis}).
\end{proof}

\begin{prop}
\label{prop:C-vs}
Additively, $A_n$ is a $\C_n$-vector space, where the $\C_n$-action on $A_n$
is given by left multiplication.
\end{prop}

\begin{proof}
The only thing to check is that if $\alpha$ and $\beta$ belong to $\C_n$
and $x$ belongs to $A_n$, then $\alpha (\beta x)$ equals 
$(\alpha \beta) x$.  This follows immediately from Lemma \ref{lem:i_n-alt}.
\end{proof}

From now on, we will often view
$A_n$ not just as an $\R$-vector space but also as a $\C_n$-vector space.
  
\begin{lemma}
\label{lem:C-linear}
Let $\phi\colon \C_{n}^\perp \map A_{n+1}$ be the map that takes $a$ to
$(a,i_n a)$.
Then one has
$\phi((p+qi_n)a)=(p+qi_{n+1})\phi(a)$ for all $p$ and $q$ in $\R$.  That is to
say, $\phi$ is complex-linear.  The same is true for the map that takes
$a$ to $(a, -i_n a)$.
\end{lemma}

\begin{proof}
Let $\alpha$ belong to $\C_{n}$, and let $a$ belong to
$\C_n^\perp$.
Write $\alpha$ as $p + qi_n$.
We want to show that $(\alpha a, i_n (\alpha a))$ equals
$(p,q)(a,i_n a)$ as elements of $A_{n+1}$.  
This follows from direct computation, using that $a^* = -a$ and
that $a$ and $i_n$ anti-commute.
\end{proof}

The element $i_n$ has some special properties not enjoyed by a typical
imaginary unit vector.

\begin{lemma}
\label{lem:i-alt}
For all $x$ in $A_n$, the associators $[x,x,i_n]$ and
$[i_n,x,x]$ both vanish.
\end{lemma}

\begin{proof}
To show that $[x,x,i_n]$ vanishes, directly compute with the inductive
definition of multiplication.  One needs to use that expressions of
the form $a + a^*$ are central because they are real.  Also, one needs
to know that $aa^*$ and $a^*a$ are equal (see Lemma~\ref{lem:norm}).

The argument for $[i_n, x, x]$ is similar.
\end{proof}

In fact, the property expressed in Lemma \ref{lem:i-alt}
determines $i_n$ uniquely, up to a sign.  Namely, if $x$ is an
imaginary unit vector in $A_n$ such that $[x,y,y] = 0$ for all $y$,
then $x = i_n$ or $x = -i_n$ \cite[Lemma 1.2]{ES}.  We also remark
that there is an automorphism of $A_n$ which fixes every element of
$A_{n-1}$ and sends $i_n$ to $-i_n$.

\begin{lemma}
\label{lem:i-quat-pair}
If $a$ is a unit vector in $\C_n^\perp$, then $a$ and $i_n$
are a quaternionic pair (see Definition \ref{defn:quat-pair}).
\end{lemma}

\begin{proof}
Since $a$ and $i_n$ are orthonormal imaginary vectors,
it suffices to show that $[i_n,a,a] = [i_n,i_n,a] = 0$.  We have
already checked these in Lemmas~\ref{lem:i_n-alt} and ~\ref{lem:i-alt}.
\end{proof}

The $\C_n$-vector space $A_n$ is not a $\C_n$-algebra.  However, we have
the following two partial results along these lines.
These lemmas are the key to computing with $i_n$.  They
allow one to do essentially any desired manipulation with
expressions involving $i_n$.

\begin{lemma}
\label{lem:i-comm}
Suppose that $x$ belongs to $\C_n^{\perp}$, and let $\alpha$ belong to $\C_n$.
For all $y$,
$(yx)\alpha = (y \alpha^*) x$ and $\alpha(xy) = x(\alpha^* y)$.
\end{lemma}

\begin{proof}
By linearity, we may assume that $\alpha$ equals $1$ or $i_n$.  The first
case is easy.  For the second case, 
compute with the inductive definition of multiplication.
\end{proof}

\begin{lemma}
\label{lem:i-comm2}
If $x$ and $y$ anti-commute and $\alpha$ belongs to $\C_n$, 
then $(\alpha x)y = -(\alpha y)x$ and $y(x \alpha) = -x(y \alpha)$.
\end{lemma}

\begin{proof}
By linearity, we may assume that $\alpha$ is either $1$ or $i_n$.
The first case is obvious.  

For the second case,
$[i_n, x+y, x+y] = 0$ by Lemma~\ref{lem:i-alt}.
Expand this by linearity to obtain
\[
[i_n, x, x] + [i_n, x, y] + [i_n, y, x] + [i_n, y, y] = 0.
\]
The first and fourth terms are zero.  Expand the other two terms to obtain
\[
(i_n x)y - i_n(xy) + (i_n y)x - i_n(yx) = 0.
\]
The second and fourth terms cancel because $xy = -yx$; the remaining two
terms give the first desired identity.

The second identity can be obtained by conjugating the first identity.
\end{proof}


\section{Hermitian inner product}
\label{sctn:herm}

This Hermitian inner product was first considered by Moreno \cite{M3}.

\begin{defn}
\label{defn:Herm}
Let $x$ and $y$ belong to $A_n$.
The \mdfn{Hermitian inner product $\langle x, y \rangle_H$} of $x$ and $y$ 
is the orthogonal projection of $xy^*$ onto $\C_n$.
\end{defn}

\begin{remark}
\label{rem:Herm}
One can check that 
$\langle x, y \rangle_H$ equals
$\langle x, y \rangle - i_n \langle i_n x, y \rangle$.  This follows
from the definition of projection once one has checked that $\langle
i_n,xy^* \rangle = -\langle i_nx,y\rangle$, and this last identity is
an easy consequence of Lemma~\ref{lem:adjoint}.
\end{remark}

\begin{prop}
\label{prop:Herm}
Definition~\ref{defn:Herm} satisfies the usual properties of a Hermitian
inner product.
\end{prop}

\begin{proof}
First, the inner product is additive in both variables.  

Second, we show that $\langle \alpha x, y \rangle_H$ equals
$\alpha \langle x, y \rangle_H$ for all $\alpha$ in $\C_n$.
By linearity, we may assume that $\alpha$ is either $1$ or $i_n$.
The formula is obvious if $\alpha$ is 1, so we may assume that
$\alpha$ equals $i_n$.  In this case, $\langle i_n x, y \rangle_H$
equals $\langle i_n x, y \rangle - i_n \langle -x, y \rangle$
by Remark~\ref{rem:Herm} and the fact that $i_n(i_n x)$ equals $-x$.
Now this expression equals 
$i_n ( \langle x,y \rangle - i_n \langle i_n x,y \rangle)$, which
is equal to $i_n \langle x, y \rangle_H$ as desired.

Next we show that $\langle x, y \rangle_H$ and $\langle y, x \rangle_H$
are conjugate.
This follows from the fact that $xy^*$ and $yx^*$ are conjugates; therefore,
their projections onto $\C_n$ are also conjugate.

Finally, we show that the inner product is positive-definite.
For any $x$, $xx^*$ is real.  Therefore, 
$\langle x, x \rangle_H$ equals $xx^*$, which equals $|| x ||^2$.
\end{proof} 

Recall from Definition~\ref{defn:L_x} that $L_x$ is the linear map
$A_n \map A_n$ given by left multiplication by $x$.

\begin{lemma}
\label{lem:antilinear}
If $x$ belongs to $\C_n^\perp$, then $L_x$ is conjugate-linear
in the sense that $L_x(y+z) = L_x(y) + L_x(z)$ and 
$L_x(\alpha y) = \alpha^* L_x(y)$ for $\alpha$ in $\C_n$.
\end{lemma}

\begin{proof}
This is simply a restatement of Lemma~\ref{lem:i-comm}.
\end{proof}

\begin{lemma}
\label{lem:L_x-Herm}
If $x$ belongs to $\C_n^\perp$, then
$L_x$ is anti-Hermitian in the sense that 
$\langle L_x y, z \rangle_H = - \langle y, L_x z \rangle_H^*$ .
\end{lemma}

\begin{proof}
Start with $\langle xy, z \rangle_H$, which equals
$\langle xy, z \rangle - i_n \langle i_n (xy), z \rangle$ by
Remark~\ref{rem:Herm}.  By 
Lemma~\ref{lem:i-comm}, this equals
$\langle xy, z \rangle + i_n \langle x (i_n y), z \rangle$.  Now
using Corollary~\ref{cor:L_x-sym}, this expression
equals
$-\langle y, xz \rangle - i_n \langle i_n y, xz \rangle$.  This is
the negative conjugate of $\langle y, xz \rangle_H$, as desired.
\end{proof}

We record the following two results about conjugate-linear anti-Hermitian
maps for later use.

\begin{lemma}
\label{lem:sing}
Suppose that $V$ is an odd-dimensional $\C$-vector space with a
nondegenerate Hermitian inner product.
If $L$ is a conjugate-linear anti-Hermitian endomorphism of $V$, then
$L$ is singular.
\end{lemma}

\begin{proof}
Choose a basis for $V$, and
identify the elements of $V$ with column vectors.  Because $L$
is conjugate-linear, there exists a
complex matrix $A$ such that $Lx = Ax^*$.  The Hermitian inner
product on $V$ is given by $\langle x, y \rangle = x^T H y^*$
for some Hermitian matrix $H$.  This means that $H^T$ equals $H^*$.
Also, because the inner product is nondegenerate, the matrix $H$ is
invertible.

The inner 
product $\langle Lx, y \rangle$ is equal to $(x^*)^T A^T H y^*$.
On the other hand, 
$-\langle x, Ly \rangle^*$ is equal to $-(x^*)^T H^* A y^*$.
Since $L$ is anti-Hermitian, these two expressions are equal
for all $x$ and $y$.
This means that $A^T H$ equals $-H^* A$, which is equal to
$-(A^T H)^T$ because $H^*$ equals $H^T$.

Write $B = A^T H$.  Then $B = -B^T$.  Since $H$ is invertible,
$A$ is singular if and only if $B$ is singular.
If $n$ is the dimension of $V$, then
\[
\det (B) = \det (-B^T) = (-1)^n \det(B^T) = (-1)^n \det(B).
\]
When $n$ is odd, this implies that $\det B$ is zero.
\end{proof}

\begin{lemma}
\label{lem:codim}
Suppose that $V$ is a $\C$-vector space with a
nondegenerate Hermitian inner product, and let 
$L$ be a conjugate-linear anti-Hermitian endomorphism of $V$.
Then the $\C$-codimension of $\ker L$ in $V$ is even.
\end{lemma}

\begin{proof}
Let $K=\ker L$.
Consider the space $W$ of all vectors $y$ such that $\langle y, z \rangle_H$
equals zero for all $z$ in $K$.  
The dimension of $W$ is the same as the codimension of $K$ because $W$
is the orthogonal complement of $K$.
Thus, we want to show that $W$ is even-dimensional.

First we will show that $L$ restricts to a map from $W$ to itself.
Suppose that $w$ belongs to $W$; then $\langle w, z \rangle_H$ equals
zero for all $z$ in $K$.  Now $\langle L w, z \rangle_H$ equals
$-\langle w, L z \rangle_H^*$.  For $z$ in $K$, $L z = 0$ by definition.
Therefore, $\langle L w, z \rangle_H$ equals zero for all $z$ in $K$.
This implies that $L w$ belongs to $W$.

Now the restriction of $L$ to $W$ is still conjugate-linear and anti-Hermitian.
Moreover, it is also non-singular because we have ensured that $W$
does not meet the kernel of $L$.  Therefore, Lemma~\ref{lem:sing} 
tells us that the dimension of $W$ must be even.
\end{proof}


\section{Automorphisms of $A_3$}
\label{sctn:Aut-A_3}

In this section we summarize some facts about automorphisms of $A_3$.
These will be used to simplify some calculations later on.  These
rather old results are due to \'{E}lie Cartan~\cite{ecartan}.  See
\cite{ES} for a generalization of this result to all $A_n.$

\medskip

Note that $A_3$ is generated as an $\R$-algebra by the three elements
$i_1$, $i_2$, and $i_3$.  In order to keep the notation more readable,
in this section we will refer to these three elements as $i$, $j$, and $t$
respectively.

To construct an algebra map $\phi:A_3 \map A_3$, one just needs to
specify $x=\phi(i)$, $y=\phi(j)$, and $z=\phi(t)$.  We know that $x$
and $y$ are a quaternionic pair.  
This implies that $x$ and $y$ are orthogonal imaginary unit vectors.
We also know that $z$ anti-commutes with $x$, $y$, and $xy$.  This means
that $z$ must be an imaginary unit vector that is orthogonal to $x$,
$y$, and $xy$.  

It turns out that these conditions on $x$, $y$, and $z$ are sufficient
to guarantee that $\phi$ is an $\R$-algebra automorphism.  The proof is a
straight-forward computation.  We summarize the last few paragraphs in
the following theorem.

\begin{thm}
\label{thm:oct-sub}
There is a one-to-one correspondence between
algebra automorphisms $A_3 \map A_3$ and ordered triples $(x,y,z)$ of
pairwise orthogonal imaginary unit vectors in $A_3$ such that $z$ is
also orthogonal to $xy$.
\end{thm}

Using this concrete description of maps from $A_3$ to $A_3$, we can
describe the automorphism group of $A_3$.
It is a 14-dimensional Lie group $G_2$ 
that belongs to a fiber bundle
\[
S^3 \map G_2 \map V_2(\R^7),
\]
where $V_2(\R^7)$ is the Stiefel manifold of ordered pairs of orthonormal
vectors in $\R^7$.

By Theorem~\ref{thm:oct-sub}, an automorphism of $A_3$ corresponds to a
triple $(x,y,z)$ of imaginary vectors of norm 1 
that are pairwise orthogonal and
such that $z$ is orthogonal to $xy$.  
The map $G_2 \map V_2(\R^7)$ takes $(x,y,z)$ to $(x,y)$.  Note that
$x$ and $y$ belong to $\R^7$ because they are imaginary.  By assumption,
they are orthogonal unit vectors.  This shows that $(x,y)$ always
belongs to $V_2(\R^7)$.

To compute the fiber of $p$, we need to describe the space of imaginary
vectors $z$ of norm $1$ that are orthogonal to $x$, $y$, and $xy$.  
The orthogonality condition leaves a $4$-dimensional subspace of possibilities
for $z$.  The norm condition guarantees that $z$ belongs to a 
$3$-dimensional sphere.

Observe, in particular, that the automorphism group of $A_3$ acts
transitively on the set of imaginary unit vectors of $A_3$.

\section{Associators}

Ultimately, we are interested in understanding the annihilators of
various vectors; i.e., given $x$, we want to describe all $y$ such
that $xy = 0$.  It turns out that in order to do this, we will need to
understand the spaces of associators and
anti-associators (see Proposition~\ref{prop:ass'(a,b)} for the connection).

\begin{defn}
\label{defn:A}
For any pair of vectors $x$ and $y$, let \mdfn{$A_{x,y}$} be the linear
endomorphism of $A_n$ given by the formula $A_{x,y}(z) = [x,z,y]$.
Also, let \mdfn{$A'_{x,y}$} be the linear endomorphism of $A_n$
given by the formula $A'_{x,y}(z) = (xz)y + x(zy)$.
\end{defn}

The $A$ stands for ``associator'', of course.
Another way to denote $A_{x,y}$ is as $R_y L_x - L_x R_y$.  Similarly,
$A'_{x,y}=R_y L_x + L_x R_y$.

\begin{lemma}
\label{lem:A}
If $x$ and $y$ are imaginary, then $A_{x,y}$ is antisymmetric (in the
sense that $\langle A_{x,y} z, w \rangle$ equals 
$-\langle z, A_{x,y} w \rangle$) and $A'_{x,y}$ is symmetric.
\end{lemma}

\begin{proof}
Use that $A_{x,y} = R_y L_x - L_x R_y$
and $A'_{x,y} = R_y L_x + L_x R_y$,
together with the
fact from Corollary~\ref{cor:L_x-sym} that $L_x$ and $R_y$
are both antisymmetric.
\end{proof}

\begin{defn}
\label{defn:ass}
Given two vectors $a$ and $b$ in $A_n$, the \mdfn{associator $\Ass[a,b]$}
is the kernel of $A_{a,b}$.
The \mdfn{anti-associator $\Ass'[a,b]$} 
is the kernel of $A'_{a,b}$.
\end{defn}

Throughout this section, we will typically assume that $a$ and $b$ are
a quaternionic pair.  Recall that this means that $a$ and $b$ are 
orthogonal imaginary unit vectors
such that the subalgebra $\HH\langle a,b \rangle$ generated by $a$ and $b$
is isomorphic to the quaternions (see Definition \ref{defn:quat-pair}).  

\begin{lemma}
\label{lem:Ass-split}
Let $a$ and $b$ be a quaternionic pair, and let $L$ belong to
the subalgebra of all $\R$-linear endomorphisms of $A_n$ that is generated
by $L_a$, $L_b$, $R_a$, and $R_b$.  Then the kernel
of $L$ splits as
\[
\ker L = \big( \ker L \cap \HH\langle a,b \rangle \big) \oplus 
\big( \ker L \cap \HH\langle a,b \rangle^\perp \big). 
\]
\end{lemma}

\begin{proof}
The desired splitting follows from the fact that $L$ restricts to
endomorphisms of $\HH\langle a,b \rangle$ and $\HH\langle a,b \rangle^\perp$.
To see why
this is true, note that $L_a$, $L_b$, $R_a$, and $R_b$ restrict
to endomorphisms of $\HH\langle a,b \rangle$; therefore, $L$ also restricts
to an endomorphism of $\HH\langle a,b \rangle$.

On the other hand, Lemma~\ref{lem:subalg} indicates that $L_a$,
$L_b$, $R_a$, and $R_b$ restrict to an endomorphism of 
$\HH\langle a,b \rangle^\perp$,
so $L$ does also.
\end{proof}

Lemma~\ref{lem:Ass-split} can be applied to the maps $A_{a,b}$ and
$A'_{a,b}$ to obtain splittings of $\Ass[a,b]$ and $\Ass'[a,b]$.  Below, we
will also need to apply it to other maps.

\begin{lemma}
\label{lem:Ass-oct}
Let $a$ and $b$ be a quaternionic pair in $A_3$ (i.e., $a$ and $b$ are any
orthonormal pair of imaginary vectors).  Then $\Ass[a,b]$ equals 
$\HH\langle a,b \rangle$
and $\Ass'[a,b]$ equals $\HH\langle a,b \rangle^\perp$.
\end{lemma}

\begin{proof}
Up to automorphism, it suffices to let $a = i$ and $b = j$, so 
$\HH\langle a,b \rangle$
equals $\HH$ and $\HH\langle a,b \rangle^\perp$ equals $t\HH$.
Since $\HH$ is an associative subalgebra, it is contained in $\Ass[i,j]$.

Now let $h$ be an arbitrary element of $\HH$, and compute
that $A_{i,j}(th) = 2t(ijh)$.  This shows that $\Ass[i,j]$ intersects
$\HH^\perp$ trivially.  With Lemma~\ref{lem:Ass-split}, this proves
the first claim.

For the second claim, $\HH$ meets $\Ass'[i,j]$ trivially since
$\HH$ is an associative subalgebra.
Now, for any $h$ in $\HH$, compute that $A'_{i,j}(th) = 0$.  This shows that
$\HH^\perp$ is contained in $\Ass'[i,j]$.  With Lemma~\ref{lem:Ass-split},
this proves the second claim.
\end{proof}

Starting with the previous lemma, 
we will compute the dimensions of various associators and anti-associators
inductively.  First we need some lengthy technical computations.

\begin{lemma}
\label{lem:Ass-compute}
If $b$ is imaginary, then
\begin{enumerate}
\item
$A_{(a,0),(b,0)} (x,y) = \big( (ax)b - a(xb), (yb)a - (ya)b \big)$.
\item
$A_{(a,0),(0,b)} (x,y) = \big( b(ya) - a(by), b(ax) - (bx)a \big)$.
\item
$A'_{(a,0),(b,0)} (x,y) = \big( (ax)b + a(xb), -(yb)a - (ya)b \big)$.
\item
$A'_{(a,0),(0,b)} (x,y) = \big( b(ya) + a(by), b(ax) + (bx)a \big)$.
\end{enumerate}
\end{lemma}

\begin{proof}
Compute using the inductive definition of multiplication.
\end{proof}

\begin{lemma}
\label{lem:Ass-(a,0),(b,0)}
Let $a$ and $b$ be a quaternionic pair in $A_n$.  Then we have
\[
\Ass[(a,0),(b,0)] = \Ass[a,b] \times 
( \Ass[a,b] \cap \HH\langle a,b \rangle^\perp ).
\]
In particular, the dimension of $\Ass[(a,0),(b,0)]$ is
$2\dim \Ass[a,b] - 4$.
\end{lemma}

\begin{proof}
According to part (1) of
Lemma~\ref{lem:Ass-compute}, we need to find all $x$ and $y$
such that $(ax)b - a(xb) = 0$ and $(yb)a - (ya)b= 0$.
The solution space of the first equation is $\Ass[a,b]$.

We now have to find the solutions to the second equation.
By Lemma~\ref{lem:Ass-split} applied to $R_a R_b - R_b R_a$,
the solution
space splits as the direct sum of its intersections with 
$\HH\langle a,b \rangle$
and $\HH\langle a,b \rangle^\perp$.
Some quaternionic arithmetic indicates that the solution space
intersects $\HH\langle a,b \rangle$ trivially.

Now we may assume that $y$ belongs to $\HH\langle a,b \rangle^\perp$.  This
implies by Lemma~\ref{lem:subalg} 
that $yb$ is also in $\HH\langle a,b \rangle^\perp$.  Using that
orthogonal imaginary vectors anti-commute, we compute that
$(yb)a - (ya)b = (ay)b - a(yb)$.  Thus, the solution space of the second
equation is the
intersection of $\Ass[a,b]$ with $\HH\langle a,b \rangle^\perp$.
This space has
dimension $\dim \Ass[a,b] - 4$ because $\HH\langle a,b \rangle$ is contained
in $\Ass[a,b]$.
\end{proof}

\begin{lemma}
\label{lem:Ass-(a,0),(0,b)}
Let $a$ and $b$ be a quaternionic pair in $A_n$.  Then we have
\[
\Ass[(a,0),(0,b)] = 
(\Ass'[a,b] \oplus \R \oplus \R a)\times
(\Ass'[a,b] \oplus \R b \oplus \R ab ).
\]
In particular, the dimension of $\Ass[(a,0),(0,b)]$ is equal
to $2\dim \Ass'[a,b] + 4$.
\end{lemma}

\begin{proof}
According to part (2) of Lemma~\ref{lem:Ass-compute}, 
we need to find all $x$ and $y$
such that $b(ya) - a(by) = 0$ and $b(ax) - (bx)a = 0$.
Let $K_1$ be the solution space of the first equation, and let $K_2$
be the solution space of the second equation.

To find $K_1$ and $K_2$, first note that they
split as the direct sums of their intersections with $\HH\langle a,b \rangle$
and $\HH\langle a,b \rangle^\perp$ because of Lemma~\ref{lem:Ass-split}
applied to 
$L_b R_a - L_a L_b$ and 
$L_b L_a - R_a L_b$.
Some quaternionic arithmetic indicates that the intersection of
$\HH\langle a,b \rangle$ with $K_1$ is the 2-dimensional subspace generated by
$b$ and $ab$.
Also, the intersection of $\HH\langle a,b \rangle$ 
with $K_2$ is the 2-dimensional
subspace generated by $1$ and $a$.

Now we may assume that $x$ and $y$ belong to $\HH\langle a,b \rangle^\perp$.  
By Lemma~\ref{lem:subalg}, $ya$, $ax$, and $bx$
are also in $\HH\langle a,b \rangle^\perp$.  Using that
orthogonal imaginary vectors anti-commute, we compute that
$b(ya) - a(by) = (ay)b + a(yb)$ and $b(ax) - (bx)a = -(ax)b -a(xb)$.

Thus, the intersection of $K_1$ with $\HH\langle a,b \rangle^\perp$ equals
$\Ass'[a,b]$ since $\Ass'[a,b]$ is contained in 
$\HH\langle a,b \rangle^\perp$.  
The same is true
for $K_2$.  It follows that $\dim \Ass[(a,0),(0,b)]$ equals
\[
2 + \dim \Ass'[a,b] + 2 + \dim \Ass'[a,b].
\]
\end{proof}

\begin{lemma}
\label{lem:Ass'-(a,0),(b,0)}
Let $a$ and $b$ be a quaternionic pair.  Then $\Ass'[(a,0),(b,0)]$ equals
\[
\Ass'[a,b] \times (\Ass'[a,b] \oplus \HH\langle a,b \rangle).
\]
In particular, the dimension of $\Ass'[(a,0),(b,0)]$ is equal
to $2\dim \Ass'[a,b] + 4$.
\end{lemma}

\begin{proof}
According to part (3) of Lemma~\ref{lem:Ass-compute}, 
we need to find all $x$ and $y$
such that $(ax)b + a(xb) = 0$ and $-(yb)a - (ya)b = 0$.
The solution space of the first equation is $\Ass'[a,b]$.

Let $K$ denote the solution space of the second equation.
By Lemma~\ref{lem:Ass-split} applied to $-R_a R_b - R_b R_a$, 
$K$
splits as the direct sum of its intersections with $\HH\langle a,b \rangle$
and $\HH\langle a,b \rangle^\perp$.
Some quaternionic arithmetic demonstrates that $\HH\langle a,b \rangle$
is contained in $K$.

Now we may assume that $y$ belongs to $\HH\langle a,b \rangle^\perp$.  This
implies that $yb$
also lies in $\HH\langle a,b \rangle^\perp$.  Using that
orthogonal imaginary vectors anti-commute, we compute that
$-(yb)a - (ya)b = (ay)b + a(yb)$.
Thus, the intersection of $K$ with $\HH\langle a,b \rangle^\perp$ equals
$\Ass'[a,b]$.  
It follows that $\dim \Ass'[(a,0),(b,0)]$ equals
\[
\dim \Ass'[a,b] + 4 + \dim \Ass'[a,b].
\]
\end{proof}

\begin{lemma}
\label{lem:Ass'-(a,0),(0,b)}
Let $a$ and $b$ be a quaternionic pair.  Then $\Ass'[(a,0),(0,b)]$ equals
\[
( (\Ass[a,b] \cap \HH\langle a,b \rangle^\perp) \oplus 
   \R b \oplus \R ab )\times
( (\Ass[a,b] \cap \HH\langle a,b \rangle^\perp) \oplus \R \oplus \R a ).
\]
In particular, the dimension of $\Ass'[(a,0),(0,b)]$ is equal
to $2\dim \Ass[a,b] - 4$.
\end{lemma}

\begin{proof}
According to part (4) of Lemma~\ref{lem:Ass-compute}, 
we need to find all $x$ and $y$
such that $b(ya) + a(by) = 0$ and $b(ax) + (bx)a = 0$.
Let $K_1$ denote the solution space of the first equation, and $K_2$
the solution space of the second equation.

To find $K_1$ and $K_2$, first note that by Lemma~\ref{lem:Ass-split}
applied to 
$L_b R_a + L_a L_b$ and 
$L_b L_a + R_a L_b$,
they split as
the direct sums of their intersections with $\HH\langle a,b \rangle$ and
$\HH\langle a,b \rangle^\perp$.  Quaternionic arithmetic easily shows that the
intersection of $\HH\langle a,b \rangle$
with $K_1$ is the 2-dimensional subspace
generated by $1$ and $a$.  Also, the intersection of 
$\HH\langle a,b \rangle$ with
$K_2$ is the 2-dimensional subspace generated by $b$ and $ab$.

Now we may assume that $x$ and $y$ belong to 
$\HH\langle a,b \rangle^\perp$.  This
implies that $ax$, $bx$, and $ya$
are also elements of $\HH\langle a,b \rangle^\perp$.  Using that
orthogonal imaginary vectors anti-commute, we compute that
$b(ya) + a(by) = (ay)b - a(yb)$ and
$b(ax) + (bx)a = -(ax)b + a(xb)$.

Thus, the intersection of $K_1$ with $\HH\langle a,b \rangle^\perp$ equals
the intersection of $\Ass[a,b]$ with $\HH\langle a,b \rangle^\perp$.
The same is true
for $K_2$.  It follows that $\dim \Ass[(a,0),(0,b)]$ equals
\[
2 + (\dim \Ass'[a,b] - 4) + 2 + (\dim \Ass'[a,b] - 4).
\]
\end{proof}

Now we are ready to construct associators and
anti-associators of various prescribed dimensions.

\begin{prop}
\label{prop:ass-exist}
Let $n \geq 3$, and 
let $d \leq 2^{n-1}$ be congruent to 4 modulo 8.
Then there exist imaginary standard basis vectors 
$a$ and $b$ in $A_n$ such that 
$\dim \Ass[a,b]$ equals $d$.  There also exist imaginary 
standard basis vectors $a'$ 
and $b'$ in $A_n$ such that $\dim \Ass'[a',b']$ equals $d$.
\end{prop}

\begin{proof}
The proof is by induction.  The base case is $n=3$, which is demonstrated
in Lemma~\ref{lem:Ass-oct}.

Any two distinct imaginary standard basis vectors $a$ and $b$ are of
course orthonormal and imaginary.  Moreover,
$[a,a,b] = [b,b,a] = 0$ by Lemma \ref{lem:alt-basis}.
Therefore, $a$ and $b$ form a quaternionic pair.

Suppose that the lemma is true for $n-1$.  
Note that $d$ equals 4 or 12 modulo 16.  
We will split the proof into these two cases.

First suppose that $d$ equals 4 modulo 16.  Write $d = 16k + 4$
with $k < 2^{n-5}$.
Choose imaginary standard basis vectors $a$ and $b$ in $A_{n-1}$ such that
$\dim \Ass[a,b]$ equals $8k + 4$; this is possible by induction
because $8k + 4 \leq 2^{n-2}$.
Note that $(a,0)$, $(b,0)$, and $(0,b)$ are imaginary standard basis
vectors of $A_n$.
Then Lemma~\ref{lem:Ass-(a,0),(b,0)} says that $\dim \Ass[(a,0),(b,0)]$ 
equals $2(8k+4) - 4$, which is the same as $d$.
Also, Lemma~\ref{lem:Ass'-(a,0),(0,b)} tells us that
$\dim \Ass'[(a,0),(0,b)] = d$.

Now suppose that $d$ equals 12 modulo 16.  Write $d = 16k + 12$
with $k < 2^{n-5}$.
Choose imaginary standard basis vectors $a$ and $b$ in $A_{n-1}$ such that
$\dim \Ass'[a,b]$ equals $8k+4$; this is possible by induction
because $8k+4 \leq 2^{n-2}$.
Again, $(a,0)$, $(b,0)$, and $(0,b)$ are imaginary standard basis vectors.
Then Lemma~\ref{lem:Ass-(a,0),(0,b)} tells us that $\dim
\Ass[(a,0),(0,b)]$ equals $2(8k+4)+4$, which is just $d$.  Also,
Lemma~\ref{lem:Ass'-(a,0),(b,0)} says that
$\dim \Ass'[(a,0),(b,0)]$ equals $d$.
\end{proof}

\section{General properties of zero-divisors}

The first few basic results of this section have been known for around
a decade to most workers in the subject, certainly including Moreno
\cite{M1} and Khalil-Yiu \cite{KY}.

\begin{defn}
\label{defn:ann}
A \mdfn{zero-divisor} is a non-zero vector $x$ such that $xy = 0$
for some non-zero vector $y$.  
The \mdfn{annihilator $\Ann(x)$} of a vector $x$ in $A_n$ is the
kernel of $L_x$.
\end{defn}

Note that $\Ann(x)$ is non-zero if and only if $x$ is zero-divisor.
It turns out (by Corollary~\ref{cor:left=right} below) that $\Ann(x)$
also equals the kernel of $R_x$, so we don't have to talk about `left'
and `right' annihilators.

\begin{lemma}
\label{lem:ann-conj}
The following equations are equivalent:
\begin{enumerate}
\item
$xy = 0$.
\item
$x^*y = 0$.
\item
$xy^* = 0$.
\end{enumerate}
\end{lemma}

\begin{proof}
For the equivalence between (1) and (2),
we need to show that if $xy = 0$, then $x^* y$ is also zero.
Compute $||x^*y||^2=\Re((x^*y)(x^*y)^*) = \Re((x^*y)(y^*x))$.
Using Lemmas~\ref{lem:Re-assoc} and~\ref{lem:Re-comm}, this is 
equal to $\Re( (x(x^*y))y^*)$.  

Now a straightforward calculation shows that $x(x^* y)$ always equals
$x^*(xy)$, which is zero because $xy = 0$.  Thus $x^*y$ is a vector
of norm zero and hence is zero.

The same argument shows that (1) and (3) are equivalent.
\end{proof}

\begin{cor}
\label{cor:left=right}
For any $x$ and $y$, we have $xy = 0$ if and only if $yx = 0$.
\end{cor}

\begin{proof}
Suppose that $xy=0$.  By conjugation, $y^* x^* = 0$.  Now Lemma
\ref{lem:ann-conj} (the equivalence of (2) and (3)) implies that $yx = 0$.
\end{proof}

\begin{lemma}
\label{lem:ann-imag}
Every zero-divisor in $A_n$ is imaginary.
\end{lemma}

\begin{proof}
Suppose that $xy = 0$ for some non-zero $x$ and $y$.
Then $x^*y$ also equals zero by Lemma~\ref{lem:ann-conj}, so
$(x+x^*)y = 0$.  Since $x+x^*$ is real and $y$ is non-zero, this
shows that $x+x^*$ is zero.  In other words, $x$ is imaginary.
\end{proof}

\begin{lemma}
\label{lem:ann-i}
Every zero-divisor in $A_n$ is orthogonal to $\C_n$.  So
for any non-zero $x$, $\Ann(x)$ is orthogonal to $\C_n$.
\end{lemma}

\begin{proof}
The second statement follows directly from the first.  
Suppose that $(a,b)$ is a zero-divisor.  
By Lemma~\ref{lem:ann-imag}, $a$ is imaginary; this means that
$(a,b)$ is orthogonal to $1$.

Now 
$ac - d^*b = 0$ and $da + bc^* = 0$ for some non-zero $(c,d)$.
Compute that $(-d^*, c)(b,a)$ equals
$(-d^*b - a^*c, -ad^* + cb^*)$.  
Using that $a^* = -a$, this equals
$(ac-d^*b, (da + bc^*)^*)$, which is zero.  Therefore, $(b,a)$ is
also a zero-divisor.  By Lemma~\ref{lem:ann-imag}, $b$ is imaginary.
Thus $(a,b)$ is orthogonal to $i_n$.
\end{proof}

\begin{lemma}
\label{lem:ann-C-vs}
For any non-zero $x$ in $A_n$, $\Ann(x)$ is a $\C_n$-vector space.
\end{lemma}

\begin{proof}
All we have to do is show that $\ker L_x$ is closed under
left multiplication by an element $\alpha$ of $\C_n$.  
By the previous
lemma, we may assume $x$ is orthogonal to $\C_n$.  

Lemma~\ref{lem:antilinear} says that $L_x$ is conjugate-linear.
If $L_x(y) = 0$, then $L_x(\alpha y) = \alpha^* L_x(y) = 0$.  Therefore,
if $y$ belongs to $\ker L_x$, then so does $\alpha y$.
\end{proof}

The previous lemma implies that the real dimension of $\Ann(x)$
is always a multiple of 2.  Soon we will show that
the real dimension of $\Ann(x)$ is in fact a multiple of 4.

\begin{lemma}
\label{lem:equal-ann}
For any non-zero $x$ in $A_n$ and any non-zero $\alpha$ in $\C_n$,
$\Ann(x)$ equals $\Ann(\alpha x)$.
\end{lemma}

\begin{proof}
We need to show that $xy = 0$ if and only if $(\alpha x)y = 0$.

Suppose that $xy = 0$.
We have that $x$ belongs to $\Ann(y)$, which means that $\alpha x$
belongs to $\Ann(y)$ by Lemma~\ref{lem:ann-C-vs}.  It follows
that $(\alpha x)y = 0$.

On the other hand, suppose that $(\alpha x)y = 0$.  
We have that $\alpha x$ belongs to $\Ann(y)$, which means that 
$x = \alpha^*\alpha x/||\alpha||^2$ 
belongs to $\Ann(y)$ by Lemma~\ref{lem:ann-C-vs}.  It follows
that $xy = 0$.
\end{proof}

The following result was originally proven by Moreno---see
\cite[Cor. 1.17]{M1}.

\begin{thm}
\label{thm:4-dim}
Let $n \geq 2$.
For any $x$ in $A_n$, the real dimension of $\Ann(x)$ is a multiple of $4$.
\end{thm}

\begin{proof}
If $x$ is zero, then $\Ann(x)$ equals $A_n$, so it has dimension $2^n$.
By the assumption on $n$, this is a multiple of $4$.  Likewise, if $x$
is not a zero-divisor, then $\Ann(x) = 0$.

If $x$ is a zero-divisor, then Lemma~\ref{lem:ann-i} says that $x$
belongs to $\C_n^\perp$.  Under these conditions, Lemmas~\ref{lem:antilinear}
and ~\ref{lem:L_x-Herm} imply that the map $L_x$ is conjugate-linear
and anti-Hermitian.  Lemma~\ref{lem:codim} implies that the complex
codimension of the kernel of $L_x$ is even.  Since the complex
dimension of $A_n$ is $2^{n-1}$, this implies that the complex
dimension of the kernel of $L_x$ is also even.

Thus, the kernel is an even-dimensional $\C_n$-vector space, so its
real dimension is a multiple of $4$.
\end{proof}

\begin{lemma}
\label{lem:zd-upbound}
Let $a$ and $b$ belong to $A_{n-1}$.  The dimension of 
$\Ann(a,b)$ is at most $2^{n-1} - 2 + \dim ( \Ann(a) \cap \Ann(b) )$.
\end{lemma}

\begin{proof}
Recall that Lemma~\ref{lem:ann-i} tells us that $\Ann(a,b)$ is a
subspace of the $(2^n-2)$-dimensional space $\C_n^\perp$. 
Let $W$ be the subspace of $\C_n^\perp$ consisting of all vectors of
the form $(c,0)$ with $c$ imaginary in $A_{n-1}$.  This is
a $(2^{n-1}-1)$-dimensional subspace of $\C_n^\perp$.

Let us investigate the intersection $\Ann(a,b) \cap W$, which
consists of vectors
$(c,0)$ such that $(a,b)(c,0) = 0$.  This means that 
$ac$ and $bc^*$ are zero.  In other words, 
$\Ann(a,b) \cap W$ is equal to $(\Ann(a) \cap \Ann(b)) \times 0$.

Now $W$ and $\Ann(a,b)$ are both subspaces of $\C_n^\perp$.  Therefore,
\[
\dim W + \dim \Ann(a,b) \leq \dim \C_n^\perp + \dim (W \cap \Ann(a,b)).
\]
Plugging in what we know, we get
\[
2^{n-1} - 1 + \dim \Ann(a,b) \leq 2^n - 2 + \dim (\Ann(a) \cap \Ann(b))
\]
Now just simplify the inequality to obtain
\[
\dim \Ann(a,b) \leq 2^{n-1} - 1 + \dim (\Ann(a) \cap \Ann(b)).
\]
Finally, $\Ann(a,b)$ and $\Ann(a) \cap \Ann(b)$ are both complex
vector spaces (see Lemma~\ref{lem:ann-C-vs}), 
so their real dimensions are even.
\end{proof}

\begin{prop}
\label{prop:zd-upbound}
Let $n \geq 2$.
For any non-zero $x$ in $A_n$, the real dimension of $\Ann(x)$ is 
at most $2^n - 4n + 4$.
\end{prop}

\begin{proof}
The proof is by induction.  The base cases are $n=2$ and $n = 3$, 
which say that $A_2$ and $A_3$ have no zero-divisors.

Assume for induction that the proposition is true for $n-1$.
Let $x = (a,b)$, where $a$ and $b$ belong to $A_{n-1}$.
By the induction assumption, we know that 
$\dim ( \Ann(a) \cap \Ann(b) ) \leq 2^{n-1} - 4(n-1) + 4$.
Now Lemma~\ref{lem:zd-upbound} implies that
$\dim \Ann(x) \leq 2^{n-1} - 2 + 2^{n-1} - 4(n-1) + 4$,
which simplifies to the inequality
$\dim \Ann(x) \leq 2^n - 4n + 6$.
Finally, recall from Theorem~\ref{thm:4-dim} that
$\dim \Ann(x)$ is a multiple of 4.
\end{proof}

The above result shows that $\Ann(x)$ has dimension at most $4$ in
$A_4$---i.e., that every zero-divisor in $A_4$ has a 4-dimensional
annihilator.  However, as $n$ increases the proposition seems to
become very weak: it gives a linear lower bound on the codimension of
each annihilator, while the dimension of $A_n$ grows exponentially.
For example,
it says that in $A_6$ the kernel of $L_x$ has dimension at most $44$.
This, in conjunction with the fairly naive method used to prove
Proposition~\ref{prop:zd-upbound}, makes it seem rather surprising
that the upper bound it establishes is in fact sharp.  
We will show  that there exists an element $x$ of $A_n$ such that
the real dimension of $\Ann(x)$ is equal to $2^n - 4n + 4$.  Moreover,
we will show in Theorem \ref{thm:ann-large-exist}
that all smaller dimensions (that are multiples of 4)
also occur.

We close this section with two simple lemmas that will be needed
later.

\begin{lemma}
\label{lem:zd-equiv}
Let $a$ and $b$ be imaginary elements of $A_{n-1}$.  
The following three subsets of $A_n$ are identical:
\begin{enumerate}[(i)]
\item
$\Ann (a,b)$. 
\item The set of all $(x,y)$ such that $x$ and $y$ are imaginary,
$ax=-yb$, and $bx=ya$.  
\item The set of all $(x,y)$ such that $x$ and $y$ are imaginary,
$ax=-yb$, and $xb=ay$.  
\end{enumerate}
\end{lemma}   

\begin{proof}
The equivalence of the first two subsets comes from just writing
out $(a,b)(x,y) = (0,0)$ with the definition of multiplication.
Also, recall from Lemma~\ref{lem:ann-i} that $x$ and $y$ have to be
imaginary, so $x^* = -x$ and $y^* = -y$.

For the third set, conjugate the equation $bx=ya$, using that
$a$, $b$, $x$, and $y$ are all imaginary.
\end{proof}

\begin{lemma}
\label{lem:0-subalg}
Suppose that $B$ is a subalgebra of $A_n$ containing a vector $x$.
Then $\Ann(x)$ decomposes as
\[
(\Ann(x) \cap B) \oplus (\Ann(x) \cap B^{\perp}).
\]
\end{lemma}

\begin{proof}
Let $y$ belong to $\Ann(x)$, and write $y = y_1 + y_2$, where
$y_1$ belongs to $B$ and $y_2$ belongs to $B^\perp$.  All we have to
do is show that $y_1$ and $y_2$ also belong to $\Ann(x)$.

Since $B$ is a subalgebra, $xy_1$ belongs to $B$.
Recall from Lemma~\ref{lem:subalg} that $xy_2$ belongs to $B^{\perp}$.
Now $xy_1$ and $xy_2$ are orthogonal vectors whose sum is zero,
so they must both be zero.
\end{proof}


\section{Constructions of zero-divisors}

We now begin the task of producing
zero-divisors whose annihilators have various dimensions.

\begin{thm}
\label{thm:C-ann}
Let $a$ be a vector in $\C_n^\perp$, and
let $\alpha$ and $\beta$ be elements of $\C_n$
such that $\alpha^2 + \beta^2$ is not zero.
Then $\Ann(\alpha a, \beta a)$
is equal to $\Ann(a) \times \Ann(a)$; in particular, the dimension of
$\Ann(\alpha a, \beta a)$ is $2 \dim \Ann(a)$.
\end{thm}

\begin{proof}
First note that $\Ann(a) \times \Ann(a)$ is contained in $\Ann(\alpha
a,\beta a)$, using that $\Ann(a)\subseteq \Ann(\alpha a)$ (and
similarly for $\beta$), with
equality unless $\alpha=0$.  So we will prove the subset in the other
direction. 

Without loss of generality, we may rescale $a$ and assume that it is
a unit vector.
Recall that $a$ and $i_n$ are a quaternionic pair; they generate a subalgebra
$\HH\langle a,i_n \rangle$
of $A_n$ isomorphic to the quaternions with additive basis
consisting of $1$, $a$, $i_n$, and $i_n a$.

Consider the subalgebra 
$B = \HH\langle a,i_n \rangle \times \HH\langle a,i_n \rangle$ of $A_{n+1}$.
Note that $B$ is isomorphic to the octonions, and $(\alpha a, \beta a)$
belongs to $B$.  Since the octonions have no zero-divisors,
Lemma~\ref{lem:0-subalg} implies that $\Ann(\alpha a, \beta a)$
is contained in $B^{\perp}$.  

Our goal is to find all $x$ and $y$ satisfying the equation 
$(\alpha a, \beta a)(x,y) = (0,0)$.  
The previous paragraph says that $x$ and $y$ must belong to 
$\HH\langle a,i_n \rangle^\perp$.  In particular, $x$ and $y$ are orthogonal
to $\alpha a$ and $\beta a$.

The equation $(\alpha a, \beta a)(x,y) = (0,0)$ is equivalent to the pair of 
equations $(\alpha a)x + y (\beta a) = 0$ and $y(\alpha a) - (\beta a)
x = 0$, since $a$, $x$, and $y$ are all imaginary.

We'll work with the first equation first.  Since $\alpha a$ and $x$
anti-commute, we get
$-x(\alpha a) + y (\beta a) = 0$.
Lemma~\ref{lem:i-comm} implies that $\alpha a$ equals $a \alpha^*$
and $\beta a$ equals $a \beta^*$, so we obtain
$-x(a \alpha^*) + y (a \beta^*) = 0$.
Next use Lemma~\ref{lem:i-comm2} to get
$a(x \alpha^*) - a (y \beta^*) = 0$.
Finally, use Lemma~\ref{lem:i-comm} again and factor to obtain
$a(\alpha x - \beta y) = 0$.  Therefore, the first equation is equivalent
to the condition that $\alpha x - \beta y$ belongs to $\Ann(a)$.

For the second equation, use similar arguments to get
$-a ( \beta x + \alpha y ) = 0$.  Therefore, the second equation 
is equivalent to the condition that $\beta x + \alpha y$ belongs to $\Ann(a)$.

Since $\alpha^2 + \beta^2$ is not zero, it follows that
$x$ and $y$ both belong to $\Ann(a)$, as desired.
\end{proof}

The previous theorem handles a large class of zero-divisors of
the form $(\alpha a, \beta a)$, where $\alpha$ and $\beta$ belong to
$\C_n$ while $a$ belongs to $\C_n^\perp$.  However, it does not include
the situation where $\alpha^2 + \beta^2 = 0$, i.e., when 
$\beta$ equals $i_n \alpha$ or $-i_n \alpha$.  The following theorem
takes care of these remaining cases.

\begin{thm}
\label{thm:C-ann2}
If $a$ is a vector in $\C_n^\perp$, then $\Ann(a,i_n a)$ is equal to
\[
\{ (x,i_n x) : x \in \Ann(a) \} \oplus
\{ (y, -i_n y) : y \in \HH\langle a,i_n \rangle^\perp  \},
\]
and $\Ann(a, -i_n a)$ is equal to
\[
\{ (x,-i_n x) : x \in \Ann(a) \} \oplus
\{ (y, i_n y) : y \in \HH\langle a,i_n \rangle^\perp  \},
\]
In particular, the dimensions of $\Ann(a,i_n a)$ and of $\Ann(a,-i_n a)$
are both equal to $2^n - 4 + \dim \Ann(a)$.
\end{thm}

\begin{proof}
We prove the theorem for $\Ann(a,i_n a)$; the proof for 
$\Ann(a,-i_n a)$ is identical (or one can use the automorphism of
$A_n$ which fixes $A_{n-1}$ pointwise and interchanges $i_n$ and
$-i_n$). 

By direct computation, one can verify that $(x,i_n x)$ belongs to
$\Ann(a,i_n a)$ when $ax = 0$ and that $(y,-i_n y)$ belongs to
$\Ann(a,i_n a)$ when $y$ belongs to $\HH\langle a,i_n \rangle^\perp$.

Now suppose that $(z,w)$ belongs to $\Ann(a,i_n a)$.
Write $(z,w)$ in the form $(x,i_n x) + (y, -i_n y)$, where 
$x$ equals $(z - i_n w)/2$ and $y$ equals $(z + i_n w)/2$.
We want to show that $x$ belongs to $\Ann(a)$ and that 
$y$ belongs to $\HH\langle a,i_n \rangle^\perp$.

As in the proof of 
Theorem~\ref{thm:C-ann}, we know that $z$ and $w$ must belong 
to $\HH\langle a,i_n \rangle^\perp$.  
It follows from Lemma~\ref{lem:subalg} that
$i_n w$ is also in $\HH\langle a,i_n \rangle^\perp$.  This shows that $y$
belongs to $\HH\langle a,i_n \rangle^\perp$.

We also know from the proof of Theorem~\ref{thm:C-ann} that
$z - i_n w$ belongs to $\Ann(a)$.  That is, $x\in \Ann(a)$.
\end{proof}


\section{Anti-associators and zero-divisors}

In this section, our goal is to describe $\Ann (a,b)$ when $a$ and $b$
are a quaternionic pair of alternative vectors.  Our description will
be in terms of the anti-associator $\Ass'[a,b]$.  There appears to be
a connection between this result and some statements in~\cite{M2};
Moreno was the first person to study zero-divisors which are pairs of
alternative vectors.

\begin{prop}
\label{prop:ass'(a,b)}
Let $a$ and $b$ be a quaternionic pair of alternative vectors.  Then
$\Ann (a,b)$ is equal to
\[
\{ (x, (ax)b) : x \in \Ass'[a,b]  \}.
\]
In particular, the dimension of $\Ann(a,b)$ equals
$\dim \Ass'[a,b]$.
\end{prop}

\begin{proof}
Suppose that $(a,b)(x,y) = (0,0)$.  Since we know $x$ and $y$ must be
imaginary, this equation is equivalent to the
two equations $ax = -yb$ and $xb = ay$ (see Lemma~\ref{lem:zd-equiv}).
Multiply the first equation by $b$ on the right to obtain
$y = (ax)b$, and multiply the second equation by $a$ on the left
to obtain $y = -a(xb)$.
\end{proof}

\begin{cor}
\label{cor:ann-A_4}
Let $a$ and $b$ be orthogonal imaginary vectors in $A_3$ such that
$||a|| = ||b|| \neq 0$.
Then we have
\[
\Ann(a,b) = \{ (x, -(ab)x/||ab||) : x \in \HH\langle a,b \rangle^\perp  \}.
\]
\end{cor}

\begin{proof}
By rescaling $(a,b)$, we may assume that $a$ and $b$ are unit vectors.

Note that $a$ and $b$ are automatically alternative because every
vector in $A_3$ is alternative.  So $a$ and $b$ are 
automatically a quaternionic pair.  Therefore, 
Proposition~\ref{prop:ass'(a,b)}
applies.  We have replaced $\Ass'[a,b]$ with $\HH\langle a, b \rangle^\perp$
with the help of Lemma \ref{lem:Ass-oct}.

Now we just have to do some octonionic arithmetic
and observe that $(ax)b$ equals $-(ab)x$; here we need that $x$
is orthogonal to $\HH\langle a,b \rangle$, 
and we are using Lemma~\ref{lem:Ass-oct}.
\end{proof}


\section{Zero-divisors in $A_4$}
\label{sctn:zd-A_4}

The algebras $A_0$, $A_1$, $A_2$, and $A_3$ are all normed algebras.
Therefore, they have no zero-divisors.  However,
$A_n$ does have zero-divisors when $n\geq 4$.
The purpose of this section is to thoroughly describe the pairs
of non-zero vectors $(x,y)$ in $A_4$ such that $xy = 0$.

The results of this have been known for nearly two decades, but their
provenance is a bit complicated.  The main ingredients can be found
in~\cite{ES}, but it does not seem that Eakin and Sathaye were aware
of this.  Cohen states the main result without proof in~\cite{Co}, and
also asserts that Paul Yiu had told him about a different, also
unpublished proof.  To our knowledge, the first complete published
proofs of these results are~\cite[Corollary 2.14]{M1} and
\cite[Theorem 3.2.3]{KY}. 

\begin{prop}
\label{prop:zd-A_4}
A vector $(a,b)$ in $A_4$ is a zero-divisor if and only if $a$ and $b$ are 
orthogonal imaginary vectors such that $||a|| = ||b||$.
\end{prop}

\begin{proof}
One direction is Corollary~\ref{cor:ann-A_4}.  For the other direction,
suppose that $(a,b)$ is a zero-divisor.

First of all, Lemma~\ref{lem:ann-i} says that $a$ and $b$ are both imaginary.
Moreover, $a$ and $b$ must both be non-zero.  For example, 
Theorem~\ref{thm:C-ann} says that the dimension of $\Ann(a,0)$
equals $2\dim \Ann(a)$.  But $\Ann(a)$ is trivial because $A_3$ has
no zero-divisors, so $\Ann(a,0)$ is also trivial.  The same argument
applies to $\Ann(0,b)$.

There exist $x$ and $y$ such that $(a,b)(x,y) = (0,0)$.
Lemma~\ref{lem:zd-equiv} says that $ax = -yb$ and $bx = ya$.
Now $A_3$ is a normed algebra, so $||a||\cdot ||x|| = ||y||\cdot ||b||$
and $||b||\cdot ||x|| = ||y||\cdot ||a||$.  
Using that $a$, $b$, $x$, and $y$ are all non-zero (see the previous
paragraph), it follows that $||a|| = ||b||$ (and also
$||x|| = ||y||$).

It remains to show that $a$ and $b$ are orthogonal.
If we take the equation $ax = -yb$, multiply by $x$ on the right,
and use that $A_3$ is alternative, we obtain $||x||^2 a = (yb)x$.
Similarly, if we start with $bx = ya$,
we obtain $||y||^2 a = -y(bx)$.
In particular, $2||x||^2a$ equals $A_{y,x} b$ (using that $||x||=||y||$).

This now allows us to compute:
\[
\langle a,b\rangle = \frac{1}{2||x||^2}\langle A_{y,x}b,b\rangle = 0
\]
where the last equality follows from the anti-symmetry of $A_{y,x}$
(see Lemma~\ref{lem:A}).  Thus, $a$ and $b$ are orthogonal.
\end{proof}

Proposition~\ref{prop:zd-A_4} allows us to describe geometrically 
the space of all unit zero-divisors in $A_4$.  It is homeomorphic
to the Stiefel manifold $V_2(\R^7)$ of orthonormal pairs of vectors
in $\R^7$.  Here, $\R^7$ arises as the space of imaginary
vectors in $A_3$.

One can see directly from Proposition~\ref{prop:zd-A_4} and
Corollary~\ref{cor:ann-A_4} that
if $x$ is any vector in $A_4$, then $\Ann(x)$ has real dimension 0 or 4,
depending on whether $x$ is a zero-divisor or not.  Compare this
observation with Theorem~\ref{thm:4-dim} and Proposition~\ref{prop:zd-upbound}
above.

Corollary~\ref{cor:ann-A_4} allows one to describe geometrically
the space of all pairs of unit vectors in $A_4$ whose product is zero.
It turns out to be homeomorphic to the 14-dimensional Lie group $G_2$
(see Section~\ref{sctn:Aut-A_3}).

\section{Existence of annihilators with various dimensions}

In these last three sections we finally prove the main results stated
in the introduction.

\begin{prop}
\label{prop:ann-small}
Let $d$ be any non-negative integer less than $2^{n-1}$ that is congruent to
$0$ modulo $4$.  
There exists a vector $x$ in $A_n$ such that the dimension of $\Ann(x)$
is equal to $d$.
\end{prop}

\begin{proof}
The proof is by induction.  The base cases $n \leq 3$ are trivial.

First suppose that $d$ is congruent to 0 modulo 8.  Write $d = 8k$.
By induction, we may find a vector $y$ in $A_{n-1}$ such that
the dimension of $\Ann(y)$ is $4k$.  Then the dimension of
$\Ann(y,0)$ is $8k$ by Theorem~\ref{thm:C-ann}.

Now suppose that $d$ is congruent to $4$ modulo $8$.  
By Proposition \ref{prop:ass-exist}, 
we may choose imaginary standard basis vectors $a$ and $b$
such that $\dim \Ass'[a,b] = d$.  Then Proposition~\ref{prop:ass'(a,b)}
tells us that $\dim \Ann(a,b) = d$.
\end{proof}

\begin{thm}
\label{thm:ann-large-exist}
Let $n \geq 1$.
There exists a vector $x$ of $A_n$ such that $\Ann(x)$ has
dimension $d$ if and only if 
$0 \leq d \leq 2^n - 4n + 4$ and $d$ is congruent to
0 modulo 4.  
\end{thm}

\begin{proof}
One direction is a combination of Theorem~\ref{thm:4-dim} and
Proposition~\ref{prop:zd-upbound}.

For the other direction, the proof is by induction on $n$.  The bases cases
are $n \leq  3$, which require nothing.
Suppose the result has been proved for $A_{n-1}$, where $n\geq 3$, and
let $d$ satisfy the given conditions.  If $d < 2^{n-1}$, then Proposition
\ref{prop:ann-small} implies the existence of the desired $x$.

Now assume that $d \geq 2^{n-1}$.  Use induction 
to choose an $a$ in $A_{n-1}$ such that $\dim \Ann(a) = d - 2^{n-1} + 4$.
Note that $d - 2^{n-1} + 4 \leq 2^{n-1} - 4(n-1) + 4$ because
$d \leq 2^n - 4n + 4$.

Let $x = (a, i_{n-1} a)$.  Then Theorem~\ref{thm:C-ann} implies that
$\dim \Ann(x) = d$.
\end{proof}


\section{Top-dimensional annihilators in $A_5$}

The first few Cayley-Dickson algebras have no non-trivial
zero-divisors.  The fourth one has zero-divisors, but these have been
well-understood for some years now---largely because they are
homogeneous in a variety of ways.  One consequence of this homogeneity
is the fact that each zero-divisor has a 4-dimensional annihilator.

We have now demonstrated that no analogous fact holds for $A_n$ with
$n\geq 5$.  Indeed, we have shown exactly what dimensions of
annihilators occur in $A_n$ for all $n$; our results tell us 
that as $n$ increases, the number of possibilities for the dimension
of an annihilator in $A_n$ grows exponentially.

This indicates that the analysis of the space of zero-divisors in
$A_n$ will be quite complicated, but it also gives us a hint as to how
that analysis might be carried out.  Namely, write $\ZD(A_n)=\{x\in
A_n : ||x||=1, \ \Ann(x)\neq 0\}$.  We can partition this space
into the subsets
$$\ZD_k(A_n) = \{x\in A_n: ||x||=1, \ \dim\Ann(x) = k\}$$
where $k = 0,4,\dots,2^n - 4n + 4$.  This decomposition
of $A_n$ is a stratification in the sense that
$$
\overline{\ZD_k(A_n)} = \bigcup_{k'\geq k}\ZD_{k'}(A_n),$$
where the union is disjoint.
At present, it seems that the most accessible approach to the study of
the zero-divisor locus in $A_n$ is to analyze one $\ZD_k(A_n)$ at a
time.  We conclude this article with the beginning of this program,
namely a complete determination of $\ZD_{2^n - 4n + 4}(A_n)$ for all $n$.

\begin{defn}
\label{defn:top}
Let $n \geq 2$.
An element of $A_n$ is a \mdfn{top-dimensional zero-divisor} if its
annihilator has dimension $2^n - 4n + 4$.  Let $T_n$ be the space
of top-dimensional zero-divisors in $A_n$ that have norm 1.
\end{defn}

Notice that $T_n$ is nothing other than $\ZD_{2^n - 4n + 4}(A_n)$.
Proposition~\ref{prop:zd-upbound} tells us that annihilators of
zero-divisors have dimension at most $2^n - 4n + 4$, so our terminology
makes sense.  

Note that $T_2$ and $T_3$ are homeomorphic to $S^3$ and $S^7$
respectively because every annihilator in $A_2$ or $A_3$ is
zero-dimensional.  We explained in Section~\ref{sctn:zd-A_4} that
$T_4$ is homeomorphic to the Stiefel manifold $V_2(\R^7)$ of
orthonormal pairs of imaginary unit vectors in $A_3$.

In this section we will study $T_5$, which is  the space of
(necessarily imaginary)
unit vectors $x$ in $A_5$ such that $\Ann(x)$ is 16-dimensional.  This
will serve as the base case of an induction carried out in the next
section, where we describe $T_n$ for all $n\geq 5$.

\begin{lemma}
\label{lem:ann-intersect}
Let $a = (a_1, a_2)$ and $b = (b_1, b_2)$ be zero-divisors in $A_4$.
Then $\Ann(a)$ and $\Ann(b)$ intersect non-trivially if and only if 
$a_1 a_2/||a_1 a_2|| = b_1 b_2/||b_1 b_2||$.
\end{lemma}

\begin{proof}
Note that by Proposition~\ref{prop:zd-A_4} we know 
$a_1$ and $a_2$ are orthogonal and imaginary, as are $b_1$
and $b_2$.  
From Corollary~\ref{cor:ann-A_4},
the annihilators intersect if
and only if there exists a non-zero $x$ in 
$\HH\langle a_1, a_2 \rangle^\perp \cap \HH\langle b_1, b_2\rangle^\perp$
such that
$-(a_1 a_2)x/||a_1 a_2|| = -(b_1 b_2)x/||b_1 b_2||$.
Thus, we have
\[
\left( a_1 a_2/||a_1 a_2|| - b_1 b_2/||b_1 b_2|| \right) x = 0.
\]
Since the octonions have cancellation,
this equation has a non-zero solution in $x$ if and only if the left-hand
factor is zero.
\end{proof}

Recall from Proposition~\ref{prop:zd-A_4} 
that if $(a,b)$ is a zero-divisor in $A_4$, then 
$a/||a||$ and $b/||b||$ form a quaternionic pair in $A_3$.

\begin{lemma}
\label{lem:ann-equal}
Let $a = (a_1, a_2)$ and $b = (b_1, b_2)$ be zero-divisors
in $A_4$.  Then $\Ann(a)$ equals $\Ann(b)$ if and only if $a_1
a_2/||a_1 a_2|| = b_1 b_2/||b_1 b_2||$ and
$\HH\langle a_1/||a_1||,a_2/||a_2|| \rangle = 
\HH\langle b_1/||b_1||, b_2/||b_2|| \rangle$.
\end{lemma}

\begin{proof}
First note that we are free to multiply $a$ and $b$ by real scalars.
Since we already know $||a_1||=||a_2||$ and $||b_1||=||b_2||$
(by Proposition~\ref{prop:zd-A_4}),
we can assume that $||a_1|| = ||a_2|| = ||b_1|| = ||b_2|| = 1$.
This will simplify the notation somewhat.  Under this assumption, we
must show that $\Ann(a) = \Ann(b)$ if and only if $a_1 a_2 = b_1 b_2$
and $\HH\langle a_1,a_2 \rangle  = \HH \langle b_1,b_2\rangle$.

First suppose that $\Ann(a)$ and $\Ann(b)$ are equal.
It follows from Corollary~\ref{cor:ann-A_4}
that $\HH\langle a_1,a_2\rangle^\perp = \HH\langle b_1,b_2\rangle^\perp$ and 
therefore $\HH\langle a_1,a_2\rangle = \HH\langle b_1,b_2\rangle$.
Also, Lemma~\ref{lem:ann-intersect} says that 
$a_1 a_2 = b_1 b_2$.

Now suppose that $a_1 a_2 = b_1 b_2$ and $\HH\langle a_1,a_2 \rangle =
\HH\langle b_1,b_2 \rangle$.  It then follows from
Corollary~\ref{cor:ann-A_4} that $\Ann(a)$ and $\Ann(b)$ are equal.
\end{proof}

\begin{prop}
\label{prop:A_5-16}
Suppose that $(a,b)$ is an element of $A_5$ such that $a$ and $b$ 
are zero-divisors in $A_4$ with $\Ann(a) = \Ann(b)$.  Then $b=\alpha
a$ for some $\alpha\in \C_4$.  One has
that $\dim \Ann(a,b)=16$ if and only if $\alpha=\pm i_4$, and $\dim
\Ann(a,b)=8$ otherwise.
\end{prop}

\begin{proof}
Let $a=(a_1,a_2)$.  After rescaling $a$ and $b$, we may assume that
$||a_1|| = ||a_2|| = 1$.  Up to an automorphism of $A_3$, we may
additionally assume that $a_1 = i$ and $a_2 = j$.

Lemma~\ref{lem:ann-equal} implies that $b_1$ and $b_2$ belong to $\HH$
and that $b_1 b_2 = || b_1 b_2 || k$.  Since $b_1$ and $b_2$ must
both be imaginary, as well as orthogonal, it follows that $b_1 = Pi - Qj$ and
$b_2 = Qi + Pj$ for some real numbers $P$ and $Q$.

Note that $i_4 a$ equals $(-j, i)$.  Therefore, $b$ equals $Pa + Qi_4
a=(P+Qi_4)a$.  Now apply
Theorem~\ref{thm:C-ann} 
to find that $\Ann
(a,(P+Qi_4)a)$ is equal to $\Ann a\times \Ann a$ (which has dimension
$8$) if $P+Qi_4\neq \pm
i_4$.  In case $P+Qi_4=\pm i_4$, Theorem~\ref{thm:C-ann2} applies.
\end{proof}

\begin{lemma}\label{lem:t5intersect}
Suppose $(a,b)$ is an element of $T_5$.  Then $a$ and $b$ are
zero-divisors in $A_4$ and $\Ann(a)$ and $\Ann(b)$ intersect
non-trivially.
\end{lemma}

\begin{proof}
By Lemma \ref{lem:zd-upbound}, 
\[
16=\dim \Ann(a,b) \leq 16 - 2 + \dim ( \Ann(a) \cap \Ann(b) ).
\]
From this it follows that $\Ann(a) \cap \Ann(b)$ is non-zero.
\end{proof}

\begin{thm}\label{thm:mainA_5}
Let $a$ and $b$ be elements of $A_4$ such that 
$\Ann a\cap \Ann b\neq 0$ and $(a,b)\neq (0,0)$.  Then 
\[
\dim \Ann(a,b)=
\begin{cases}
16 & \text{if $b=\pm i_4a$,} \\
12 & \text{if $a$ is orthogonal to $b$, $||a||=||b||$, and $b\neq \pm
i_4 a$}, \\
8 & \text{otherwise}.
\end{cases}
\]
\end{thm}

Before proving this theorem we note the following immediate corollary:

\begin{cor}\label{cor:t5}
The space $T_5$ is homeomorphic to $V_2(\R^7) \amalg V_2(\R^7)$.
\end{cor}

\begin{proof}
We will show that $T_5$ is the disjoint union of the
spaces
\[
X_+ = \{ (a,i_4 a): ||a|| = \tfrac{1}{\sqrt{2}},
\ a \ \text{is a zero-divisor in $A_4$}\}
\]
and
\[
X_- = \{ (a,-i_4 a): ||a|| = \tfrac{1}{\sqrt{2}},
a \ \text{is a zero-divisor in $A_4$}\}.
\]
As explained in Section \ref{sctn:zd-A_4}, each of these spaces is
homeomorphic to $V_2(\R^7)$.

First observe that both $X_+$ and $X_-$ are contained in $T_5$
because of Theorem \ref{thm:C-ann2}.
Next we will show that $X_+$ and $X_-$ are disjoint.  If
$(a,i_4 a) = (b, -i_4 b)$, then it follows that $2i_4 a = 0$.
Since $i_4$ is alternative, this implies that $a = 0$, which prohibits
$(a,i_4a)$ from belonging to $X_+$.

Finally, we must show that every element of $T_5$ is contained
in $X_+$ or $X_-$.
Suppose that $(a,b)$ is an arbitrary element of $T_5$.  By
Lemma~\ref{lem:t5intersect}, $a$ and $b$ are zero-divisors in $A_4$
whose annihilators intersect nontrivially.  
Since $\dim \Ann(a,b)=16$, we have by Theorem~\ref{thm:mainA_5} that
$b=\pm i_4 a$.  
\end{proof}

Our final task in this section is the following:

\begin{proof}[Proof of Theorem~\ref{thm:mainA_5}]
Let $a$ and $b$ be elements $A_4$ which are not both zero, and
whose annihilators intersect nontrivially.
If either $a$ or $b$ is zero then we are in the third case from the
statement of the theorem, and the fact that the annihilator is 
$8$-dimensional follows from Theorem~\ref{thm:C-ann}.
So we may as well assume both $a$ and $b$ are nonzero.  

We can write $a = (a_1,a_2)$ and $b = (b_1,b_2)$, where $a_1,a_2,b_1$,
and $b_2$ are all octonions.
Since $a$ and $b$ are zero-divisors, 
Proposition \ref{prop:zd-A_4} implies 
that $a_1$, $a_2$, $b_1$, and $b_2$ are all imaginary, that
$a_1$ and $a_2$ are perpendicular, that $b_1$ and $b_2$ are
perpendicular, and that $||a_1|| = ||a_2||$ and $||b_1|| = ||b_2||$.
After scaling $a$ and $b$ by the same constant, we can assume
$||a_1||=||a_2||=1$.  
Since $a_1$ and $a_2$ are orthogonal imaginary unit vectors, 
we know that $a_1 a_2$ is also an imaginary unit vector.
Therefore, up to automorphism of $A_3$, we may assume that
$a_1 a_2 = k$, which implies by Lemma~\ref{lem:ann-intersect} that
$b_1 b_2$ is a scalar multiple of $k$.
Note that $a_1$, $a_2$, $b_1$, and $b_2$ are all orthogonal to $k$
since $a_1$ and $a_2$ are orthogonal to $a_1 a_2$, while
$b_1$ and $b_2$ are orthogonal to $b_1 b_2$.

Consider the linear span $V$ of $\{1,a_1,a_2,b_1,b_2,a_1 a_2,b_1 b_2\}$;
its dimension is at most 6.
Thus, up to automorphism, we may assume that
$i$ is orthogonal to $V$.
Since $V$ is invariant under left
multiplication by $k$, 
$j$ is orthogonal to $V$ as well.
We have now shown that $a_1$, $a_2$, $b_1$, and $b_2$ are all orthogonal
to $\HH$.  Thus, once again up to automorphism, we may assume that
$a_1 = t$.  This implies $a_2 = kt$, as $a_1a_2=k$.  

Since $b_1$ is orthogonal to $\HH$, 
we have that $b_1 = \alpha t$ for some element
$\alpha$ of $\HH$.  It then follows that $b_2 = (\alpha k)t$
because $b_1 b_2$ is a scalar multiple of $k$ and because $||b_1|| = ||b_2||$.
So we have
$a=(t,kt)$ and $b=(\alpha t,(\alpha k)t)$.  To relate our present situation
to the three cases in
the statement of the theorem, note that $b=\pm i_4a$ if and only if
$\alpha=\pm k$; $a$ is orthogonal to $b$ if and only if $\alpha$ is
imaginary; and $||a||=||b||$ if and only if $||\alpha||=1$.

Consider now the subalgebra $\HH$ of $A_3;$ we may use it to form the algebra
$$\Omega = \HH\times\HH\subset A_4.
$$
This in turn gives rise to a 16-dimensional subalgebra
$\Omega\times\Omega$ of $A_5$.  By assumption, $(a,b)$ is an element
of $\Omega^\perp\times\Omega^\perp$.  Moreover, one easily sees that in
addition to the automatic
$$(\Omega\times\Omega)\cdot(\Omega^\perp\times\Omega^\perp)\subset
\Omega^\perp\times\Omega^\perp,
$$ 
we have the formula
$$(\Omega^\perp\times\Omega^\perp)\cdot(\Omega^\perp\times\Omega^\perp)\subset
\Omega\times\Omega.
$$
Similarly to Lemma \ref{lem:0-subalg},
this in turn implies that $\Ann(a,b)$ is isomorphic to
$$(\Ann(a,b)\cap(\Omega\times\Omega)) \oplus
(\Ann(a,b)\cap(\Omega^\perp\times\Omega^\perp)).
$$

In other words, we must solve the two equations
\begin{myequation}\label{eq:omegasols}
(t,kt,\alpha t,(\alpha k)t)\cdot(x,y,z,w) = (0,0,0,0)
\end{myequation}
and
\begin{myequation}\label{eq:omegaperpsols}
(t,kt,\alpha t,(\alpha k)t)\cdot(xt,yt,zt,wt) = (0,0,0,0),
\end{myequation}
where $x$, $y$, $z$, and $w$ belong to $\HH$.

At this point, the proof becomes rather unpleasant.  First of all, we
can expand out equation (\ref{eq:omegasols}). We do this by using the
inductive 
definition of multiplication in Cayley-Dickson algebras twice, and
obtain
\begin{align*}
x^* - ky^* - \alpha z^* + \alpha kw^* & =  0  \\
 	kx +  y   - \alpha kz   - \alpha w & = 0 \tag{S1}\\
 	\alpha x + \alpha ky^* +  z   + kw^* & = 0 \\
 	\alpha kx^* - \alpha y + kz^* - w & = 0.
\end{align*}
When doing the same for equation (\ref{eq:omegaperpsols}) we get
\begin{align*}
-x^*  + ky  - \alpha^*z - w^*\alpha k & = 0 \\
 	x^*k  - y  - z^*\alpha k + \alpha^*w &= 0 \tag{S2}\\
 	x^*\alpha - k\alpha^*y - z  -  w^*k &= 0 \\
 	-x^*\alpha k + \alpha^*y - z^*k + w & =0.
\end{align*}
Each of these systems corresponds to a system of $16$ real equations
in $16$ unknowns, and a little thought shows that the coefficient
matrices are negative transposes of each other (see the very end of
the proof for more information about this).  So the solution
spaces in (\ref{eq:omegasols}) and (\ref{eq:omegaperpsols}) have the
same dimension.  

We will concentrate on solving (\ref{eq:omegasols}).  We know that all
the solutions have $x$ and $z$ imaginary 
by Lemma \ref{lem:ann-i}, because $((x,y),(z,w))$
will be a zero-divisor in $A_5$. We can take advantage of this fact to
simplify the four quaternionic equations in (S1).  We obtain the new system
\begin{align*}
-x - ky^* + \alpha z + \alpha kw^* & =  0 \\
 	kx +  y   - \alpha kz   - \alpha w & = 0 \\
 	\alpha x + \alpha ky^* +  z   + kw^* & = 0 \\
 	-\alpha kx - \alpha y - kz - w & = 0 \\
        \Re(x)=\Re(z) & =0.
\end{align*}
By adding appropriate multiples of the first and third equations
(and of the second and fourth equations), the system can be simplified to
\begin{align*}
-x - ky^* + \alpha z + \alpha kw^* & =  0 \\
 	kx +  y   - \alpha kz   - \alpha w & = 0 \\
(\alpha^2+1)z + (\alpha^2+1)kw^* &=0 \\
(\alpha^2+1)kz + (\alpha^2+1)w &= 0 \\
        \Re(x)=\Re(z) & =0.
\end{align*}

\noindent
{\bf Case 1: \mdfn{$\alpha^2+1\neq 0$}}.
In this case we may cancel $\alpha^2+1$ from the last two equations
(since $\HH$ is a division algebra) and obtain $z+kw^*=0=kz+w$.
Together with $\Re(z)=0$, this is equivalent to $w\in \langle i,j\rangle$ and
$z=kw$.  Plugging this into the first two equations then gives
$x+ky^*=0=kx+y$.  Together with $\Re(x)=0$, the same analysis shows that
$y\in \langle i,j\rangle$ and $x=ky$.  So we have a $4$-dimensional
solution space for (\ref{eq:omegasols}).

\medskip
\noindent
{\bf Case 2: \mdfn{$\alpha^2+1=0$}} (equivalently, $\alpha$ is
imaginary and has norm $1$).  
In this case the third and fourth equations disappear.
We use the second equation to solve for $y$, and plug this into the
first equation.  We get (remembering that $\alpha$, $x$, and $z$ 
are imaginary):
\begin{align*}
(kxk-x) + (\alpha z+kzk\alpha) + (\alpha k w^*+kw^*\alpha) &=0 \\
\Re(x)=\Re(z) & =0,
\end{align*}
and $y$ is eliminated.  

Note that if $q$ is an imaginary quaternion of norm $1$, then
$x-qxq=2\pi_{1,q}(x)$, where $\pi_{1,q}$ denotes orthogonal
projection onto the subspace $\langle 1,q\rangle$ (it suffices
to check this claim when $q=i$).  The analysis now
divides up into two more cases.

\medskip
\noindent
{\bf Subcase 1: \mdfn{$\alpha=\pm k$}}.  The first equation becomes
$(kxk-x)\pm (kw^*k-w^*)=0$. 
So we have $\pi_{1,k}(x)=\mp\pi_{1,k}(w^*)$ and
$\Re(x)=\Re(z)=0$.  This has an $8$-dimensional solution set: $x$ and
$z$ can be any imaginary quaternions, and there are two degrees of
freedom left in choosing $w$.

\medskip
\noindent
{\bf Subcase 2: \mdfn{$\alpha \neq \pm k$}}.  Write $\alpha=r\beta+sk$
where $\beta$ is orthogonal to $k$, $||\beta||=1$, and $r,s\in \R$; so
$r^2+s^2=1$, and $r\neq 0$.  Substituting into the first equation and
re-arranging, we have
\[
(kxk-x)+ r[\beta z-(k\beta)(\beta z)(k\beta)] + 
r[\beta k w^* - \beta (\beta k w^*)\beta] + s (kw^*k-w^*)
 =0
\]
(remember that $\beta^2=-1$).  Dividing by $2$, this becomes
\[ -\pi_{1,k}(x) + r\pi_{1,\beta k}(\beta z) + r \pi_{1,\beta}(\beta k
w^*) - s \pi_{1,k}(w^*)=0.
\] 
But note that $\{1,k,\beta,\beta k\}$ is
an orthonormal basis for $\HH$, and so by separating out each
component, the above equation can be distilled into:
\begin{align*}
-\pi_1 x + r \pi_1 (\beta z) + r \pi_1 (\beta k w^*) - s \pi_1 (w^*) &= 0 \\
-\pi_k x - s\pi_k (w^*) &= 0 \\
r\pi_{\beta k} (\beta z) &= 0 \\
r\pi_\beta (\beta k w^*) &= 0.
\end{align*}
Note that 
$\pi_1 (\beta z) = \beta \pi_\beta z$, 
$\pi_1 (\beta k w^*) = \beta k \pi_{\beta k} (w^*)$,
$\pi_{\beta k} (\beta z ) = \beta \pi_k z$, and
$\pi_\beta (\beta k w^*) = \beta k \pi_k (w^*)$.
Using that $r\beta \neq 0$, together with $\Re(x)=\Re(z)=0$, we are finally
reduced to the equations
\begin{align*}
\pi_\beta z - k \pi_{\beta k} w + r^{-1}s\beta \pi_1 w &= 0 \\
-\pi_k x + s\pi_k w &= 0 \\
\pi_k z = \pi_k w = \Re(x) = \Re(z) &= 0.
\end{align*}
Since $\pi_k w = 0$, the second equation reduces to $\pi_k x = 0$.
In the end, we have three degrees of freedom for $w$, two for $x$, and
then one for $z$, yielding a six-dimensional solution space.

\medskip

We have now handled all of the cases necessary for the proof.
We will add a few comments about the two
systems (S1) and (S2).  We claimed earlier that these gave
$16\times 16$ real matrices which are negative transposes of each
other.  To see why, let $C\colon \HH \ra \HH$ denote the conjugation
operator, so $C(q)=q^*$.  We will identify $C$ with a $4\times 4$ matrix
using the standard basis for $\HH$, and where we have matrices acting
on the left.  In the same way we identify $L_q$ and $R_q$ with
$4\times 4$ matrices.

The system (S1) gives rise to a $16\times 16$ real matrix which can be
written in block form as
\[ \begin{bmatrix}
C & -L_k C & -L_\alpha C & L_\alpha L_k C \\
L_k & I & -L_\alpha L_k & -L_\alpha \\
L_\alpha & L_\alpha L_k C & I & L_k C \\
L_\alpha L_k C & -L_\alpha & L_k C & -I
\end{bmatrix},
\]
whereas the system (S2) gives the matrix
\[ \begin{bmatrix}
-C & L_k & -L_{\alpha^*} & -R_kR_\alpha C \\
R_k C & -I & -R_kR_\alpha C & L_{\alpha^*} \\
R_\alpha C & -L_kL_{\alpha^*} & -I & -R_k C\\
-R_kR_\alpha C & L_{\alpha^*} & -R_k C & I
\end{bmatrix}.
\]
Note that we have $(L_q)^T=L_{q^*}$ and $(R_q)^T=R_{q^*}$ by
Lemma~\ref{lem:adjoint}.  Also note that $CR_q=L_{q^*}C$, by the
formula $(xq)^*=q^*x^*$, and that $C$ itself is diagonal (so $C^T=C$).
Using these ideas, it follows that the two $16\times 16$ matrices are indeed
negative transposes of each other.
\end{proof}

\section{Top-dimensional zero-divisors}

The goal of this section is to completely determine the spaces $T_n$ 
for all $n$.  We have already computed some low-dimensional cases.
Our approach will be by induction, and we'll start with
several preliminary calculations.

\begin{lemma}
\label{lem:top}
Suppose that $(a,b)$ is a top-dimensional zero-divisor in $A_{n}$.  Then
$a$ and $b$ are top-dimensional zero-divisors in $A_{n-1}$, and the 
dimension of $\Ann(a) \cap \Ann(b)$ is at least $2^{n-1} - 4n + 6$.
\end{lemma}

\begin{proof}
The proof of the second claim follows from Lemma~\ref{lem:zd-upbound}
and arithmetic.

For the first claim, note that the second part implies that
the dimensions of $\Ann(a)$ and $\Ann(b)$ are at least
$2^{n-1} - 4n + 6$.  But these dimensions must be multiples of $4$
and no bigger than $2^{n-1} - 4n + 8$ by
Theorem~\ref{thm:4-dim} and Proposition~\ref{prop:zd-upbound},
so they have to be equal to $2^{n-1} - 4n + 8$.
\end{proof}

\begin{lemma}
\label{lem:not-top}
Let $a$ and $b$ be non-zero vectors in $\C_{n-2}^\perp$, and let
$x$ denote the element
$( (a, i_{n-2} a), (b, -i_{n-2} b) )$ of $A_{n}$.  The dimension of
$\Ann(x)$ is at most $2^n - 8n + 20$.
\end{lemma}

\begin{proof}
Recall that Theorem~\ref{thm:C-ann2} gives a complete description of
$\Ann(a,i_{n-2} a)$ 
and $\Ann(b, -i_{n-2} b)$.  It follows readily that 
$\Ann(a, i_{n-2} a) \cap \Ann(b, -i_{n-2} b)$ is equal to
\[
\{(x,i_{n-2}x) : x\in \Ann a\cap \HH\langle b,i_{n-2}\rangle^\perp \}
\oplus
\{(y,-i_{n-2}y) : y\in \Ann b\cap \HH\langle a,i_{n-2}\rangle^\perp
\}.
\]
So this intersection is isomorphic to
\[
( \Ann(a) \cap \HH\langle b,i_{n-2}\rangle^\perp) \oplus 
( \Ann(b) \cap \HH\langle a, i_{n-2}\rangle^\perp),
\]
which is a subspace of $\Ann(a) \oplus \Ann(b)$.
This last space has dimension at most $2^{n-1} - 8n + 24$
by Proposition~\ref{prop:zd-upbound}, so
$\Ann(a, i_{n-2} a) \cap \Ann(b, -i_{n-2} b)$ has dimension at
most $2^{n-1} - 8n + 24$.

Now Lemma~\ref{lem:zd-upbound} implies that $\Ann(x)$ has dimension at
most $2^n - 8n + 22$.  Finally, recall from Theorem~\ref{thm:4-dim} that
the dimension of $\Ann(x)$ is a multiple of 4.
\end{proof}

A similar argument shows that the dimension of the annihilator of the
element $( (a, -i_{n-2} a), (b, i_{n-2} b) )$ is also at most $2^n -
8n + 20$.  In particular, when $n \geq 5$, elements of the form $( (a,
-i_{n-2} a), (b, i_{n-2} b) )$ and $( (a, i_{n-2} a), (b, -i_{n-2} b)
)$ are never top-dimensional zero-divisors of $A_n$.

\begin{lemma}
\label{lem:H-equal}
Let $a$ and $b$ be non-zero vectors in $\C_n^\perp$.  Then 
$\HH\langle a/||a||, i_n\rangle = \HH\langle b/||b||, i_n\rangle$
if and only if $b = \alpha a$ for some element 
$\alpha$ of $\C_n$.
\end{lemma}

\begin{proof}
Without loss of generality, we may assume that $a$ and $b$ are unit vectors.

First suppose that $b = \alpha a$.  Then $b$ is an element of 
$\HH\langle a,i_n\rangle$,
so $\HH\langle b,i_n\rangle$ is contained in 
$\HH\langle a,i_n \rangle$.  On the other hand,
$a = \alpha^* b / ||\alpha||^2$, so $a$ is an element of 
$\HH\langle b,i_n\rangle$
and $\HH\langle a,i_n \rangle$ is contained in $\HH\langle b,i_n \rangle$.

Now suppose that $\HH\langle a,i_n \rangle = \HH\langle b,i_n \rangle$.
Then $b$ is a linear combination
of $1$, $i_n$, $a$, and $i_n a$.  However, since $b$ is orthogonal
to $\C_n$, $b$ is in fact a linear combination of $a$ and $i_n a$.
This is equivalent to saying that $b = \alpha a$ for some $\alpha$ in $\C_n$.
\end{proof}

\begin{lemma}
\label{lem:top-Ann-equal}
Let $a$ and $b$ be non-zero elements of $\C_n^\perp$.  The
following three conditions are equivalent:
\begin{enumerate}
\item
$\Ann(a,i_n a)$ and $\Ann(b,i_n b)$ are equal.
\item
$\Ann(a,-i_n a)$ and $\Ann(b,-i_n b)$ are equal.
\item
$b = \alpha a$ for some $\alpha$ in $\C_n$.
\end{enumerate}
\end{lemma}

\begin{proof}
First assume either condition (1) or condition (2).
By Theorem~\ref{thm:C-ann2}, we know that $\HH\langle a,i_n \rangle^\perp$
and $\HH\langle b,i_n\rangle^\perp$ are equal,
so Lemma~\ref{lem:H-equal} applies.

Now assume that $b = \alpha a$.  
By Lemma~\ref{lem:C-linear},
$(b,i_n b)$ equals $\alpha' (a, i_n a)$ for some $\alpha'$ in $\C_{n+1}$.
Similarly, $(b,-i_n b) = \alpha' (a, -i_n a)$.
Lemma~\ref{lem:equal-ann} therefore tells us that condition (1) and
condition (2) are both true.
\end{proof}

\begin{lemma}
\label{lem:codim-2}
Let $a = (a', i_{n-1} a')$ and $b = (b',i_{n-1} b')$ 
be non-zero elements of $\C_n^\perp$ such that $a'$ and $b'$
are in $\C_{n-1}^\perp$.  If 
$\Ann(a, i_n a) \cap \Ann(b,i_n b)$ has codimension at most 2
in $\Ann(a,i_n a)$, then 
$b = \alpha a$ for some $\alpha$ in $\C_n$.
\end{lemma}

A similar result holds when $a = (a',-i_{n-1} a')$ and
$b = (b', -i_{n-1} b')$.

\begin{proof}
By Theorem~\ref{thm:C-ann2}, we know that $\Ann(a, i_n a) \cap \Ann(b,i_n b)$
is isomorphic to
$( \Ann(a) \cap \Ann(b) ) \oplus 
  (\HH\langle a,i_n\rangle^\perp \cap \HH\langle b,i_n\rangle^\perp )$,
and $\Ann(a,i_n a)$ is isomorphic to 
$\Ann(a) \oplus \HH\langle a,i_n \rangle^\perp$.
Therefore, the space 
$( \Ann(a) \cap \Ann(b) ) \oplus 
  (\HH\langle a,i_n\rangle^\perp \cap \HH\langle b,i_n \rangle^\perp )$
has codimension at most 2 in 
$\Ann(a) \oplus \HH\langle a,i_n\rangle^\perp$.
All of the spaces $\Ann(a)$, $\Ann(b)$, $\HH\langle a,i_n \rangle^\perp$,
and $\HH\langle b,i_n\rangle^\perp$ are $\C_n$-vector spaces.  This means that
either $\HH\langle a,i_n\rangle ^\perp = \HH\langle b,i_n\rangle^\perp$
or $\Ann(a) = \Ann(b)$ 
(or both).  In the first case, Lemma~\ref{lem:H-equal} applies.

In the second case, Lemma~\ref{lem:top-Ann-equal} applies,
so $b' = \alpha' a'$ for some $\alpha'$ in $\C_{n-1}$.
Lemma~\ref{lem:C-linear} shows that $b = \alpha a$
for some $\alpha$ in $\C_n$.
\end{proof}

\begin{prop}
\label{prop:top}
Let $n \geq 5$.
If $x$ is a top-dimensional zero-divisor of $A_n$, then 
$x$ is of the form $(a,i_n a)$ or $(a, -i_n a)$ for some 
top-dimensional zero-divisor $a$ of $A_{n-1}$.
\end{prop}

\begin{proof}
The proof is by induction.  The base case $n = 5$ is proved in
Corollary~\ref{cor:t5}.  Notice that although the induction
establishes that $T_{n+1} = T_n \coprod T_n$ and it happens to be the
case that $T_5 = T_4\coprod T_4$, one cannot take $n=4$ as the base.
This is because the inductive step uses as a
hypothesis the way that $T_n$ is built as two copies of $T_{n-1}$, and
it is of course not the case that $T_4$ is isomorphic to a disjoint
union of two copies of $T_3$.  Aesthetically this state of affairs is
unfortunate, because it forces us to calculate $T_5$ directly in
Corollary~\ref{cor:t5}, whose proof is rather unilluminating.

Assume for induction that the proposition is true for $n-1$.  Suppose
that $x$ is a top-dimensional zero-divisor in $A_n$ and write $x =
(a,b)$ for some elements $a$ and $b$ of $A_{n-1}$. 
Lemma~\ref{lem:top} says that $a$ and $b$ are top-dimensional
zero-divisors of $A_{n-1}$, so the induction assumption says that
$a$ is of the form $(a', i_{n-2} a')$ or $(a', -i_{n-2} a')$ and
$b$ is of the form $(b', i_{n-2} b')$ or $(b', -i_{n-2} b')$.

But Lemma~\ref{lem:not-top} implies that either:
\begin{enumerate}
\item
$a = (a', i_{n-2} a')$ and $b = (b', i_{n-2} b')$, or
\item
$a = (a', -i_{n-2} a')$ and $b = (b', -i_{n-2} b')$.
\end{enumerate}
In either case, Lemma~\ref{lem:top} says that the codimension of 
$\Ann(a) \cap \Ann(b)$ in $\Ann(a)$ is at most 2,
so Lemma~\ref{lem:codim-2} applies.  It follows that
$b = \alpha a$ for some $\alpha$ in $\C_{n-1}$.

Now we know that $x$ equals $(a, \alpha a)$.  Theorem~\ref{thm:C-ann}
(and the fact that $x$ is a top-dimensional zero-divisor)
implies that $\alpha$ must be $i_{n-1}$ or $-i_{n-1}$.
Thus $x$ equals $(a, i_{n-1} a)$ or $(a, -i_{n-1} a)$.
\end{proof}

\begin{thm}
Let $n \geq 4$.  The space $T_n$ is homeomorphic to
$2^{n-4}$ disjoint copies of $V_2(\R^7)$.
\end{thm}

\begin{proof}
The proof is by induction.  The base case $n=4$ was dealt with in
Section~\ref{sctn:zd-A_4}.

Assume for induction that the theorem has been proved for $n-1$.
Proposition~\ref{prop:top} implies that $T_n$ is homeomorphic to
two disjoint copies of $T_{n-1}$.  
Note that $(a,i_{n-1} a)$ and $(b, -i_{n-1} b)$ are never equal if and $a$
and $b$ are non-zero---this follows since $i_{n-1}$ is
alternative, and therefore is not a zero-divisor.
\end{proof}

\begin{remark}
The above result prompts the following question.  Let $c$ be a
non-negative multiple of $4$, and let
$T_c(A_n)=\ZD_{2^n-4n+4-c}(A_n)$.  Does there exist a positive integer
$N$ such that for all $n> N$ the space $T_c(A_n)$ is a disjoint union
of two copies of $T_c(A_{n-1})$?  The above theorem is the case $c=0$.
\end{remark}

\bibliographystyle{amsalpha}

\begin{thebibliography}{MW}

\bibitem[Ad]{A} J. Adem, {\em Construction of some normed maps\/},
Bol. Soc. Mat. Mexicana (2) {\bf 20} (1975), no. 2, 59--75.

\bibitem[Al]{Al} A. Albert, {\em Quadratic forms permitting
composition\/}, Ann. Math. {\bf 43} (1942), no. 1, 161--177.

\bibitem[Br]{Br} R. Brown, {\em On generalized Cayley-Dickson
algebras\/}, Pacific J. Math. {\bf 20} (1967), no. 3, 415--422.


\bibitem[Ca]{ecartan}  E. Cartan, {\em Les groupes r\'{e}els simples
finis et continus\/}, Ann. Sci. \'{E}cole Norm. Sup. {\bf 31} (1914),
262--255. 

\bibitem[Co]{Co}
F. Cohen,
{\em On the Whitehead square, Cayley-Dickson algebras, and rational
functions\/}, 
Papers in honor of Jos\'e Adem,
Bol. Soc. Mat. Mexicana (2) {\bf 37} (1992), no. 1-2, 55--62.

\bibitem[Di]{Di} L. Dickson, {\em On quaternions and their
generalization and the history of the eight square theorem\/},
Ann. Math. {\bf 20} (1919), no. 3, 155--171.

\bibitem[ES]{ES}
P. Eakin and A. Sathaye, {\em On automorphisms and derivations of
Cayley-Dickson algebras\/}, J. Algebra {\bf 129} (1990), no. 2,
263--278.

\bibitem[KY]{KY}
S. Khalil and P. Yiu, {\em The Cayley-Dickson algebras, a theorem of
A. Hurwitz, and quaternions\/}, Bull. Soc. Sci. Lett. L\'od\' z S\'
er. Rech. D\' eform. {\bf 24} (1997), 117--169.

\bibitem[M1]{M1} G. Moreno, {\em The zero divisors of the Cayley-Dickson
algebras over the real numbers\/}, Bol. Soc. Mat. Mexicana (3) {\bf 4}
(1998), no. 1, 13--28.

\bibitem[M2]{M2} G. Moreno, {\em Alternative elements in the
Cayley-Dickson algebras\/}, math.RA/0404395, preprint, 2004.

\bibitem[M3]{M3} G. Moreno, {\em Hopf construction map in higher
dimensions\/}, math.AT/0404172, preprint, 2004.  

\bibitem[Sc]{Sc} R. Schafer, {\em On the algebras formed by the
Cayley-Dickson process\/}, Amer. J. Math.{\bf 76} (1954), 435--446.




\end{thebibliography}

\end{document}